\documentclass[10pt]{article}

\usepackage{graphicx}
\usepackage{graphicx,psfrag}

\def \vs {\vskip 0.3cm}

\def \vss {\vskip 0.2cm }

\def \vsss {\vskip 0.1cm }

\def \n {\noindent}

\def \E {{\bf  \, E\,}}

\def \P {{\bf P}}

\def \T {{\hbox{Tr\,}}}

\def \CV {{\cal V}}
\def \CE {{\cal E}}
\def \CC {{\cal C}}

\def \CU {{\cal U}}

\def \CT {{\cal T}}

\def \CM{{\cal M}}
\def \CW {{\cal W}}

\def \CI {{\cal I}}
\def\CL {{\cal L}}

\def \CS {{\cal S}}

\def \CF {{\cal {F}}}

\def \a {\alpha}

\def \b {\beta}

\def \g {\gamma}

\def \l {\lambda}

\def \s {\sigma}

\def \t {\tau}

\def \d {\delta}

\def \U {\Upsilon}
\def \u {\upsilon}

\def \th {\theta}

\def \vep {\varepsilon}

\def \z {\zeta}

\begin{document}

\title{  High Moments
of Large Wigner Random Matrices\\
and Asymptotic Properties of the Spectral Norm\ \footnote{ {\bf Acknowledgements:} The financial support of the research grant ANR-08-BLAN-0311-11 "Grandes Matrices 
Al\'eatoires" (France)  is gratefully acknowledged}\ 
}

\author{O. Khorunzhiy\\ Universit\'e de Versailles - Saint-Quentin, Versailles\\ FRANCE\\
{\it e-mail:} oleksiy.khorunzhiy@uvsq.fr}

\maketitle

\begin{abstract}
We consider the Wigner ensemble of  $n\times n$ 
 real symmetric random matrices $A^{(n)}$ whose entries are determined by independent identically 
distributed random variables $\{a_{ij},  i\le j\}$ that have symmetric probability distribution
with the variance $v^2$
and study the asymptotic behavior  of the spectral norm $\Vert A^{(n)}\Vert$ as \mbox{$n\to\infty$.}

We prove that if the moment $\E \vert a_{ij}\vert^{12+2\d_0}$ with any strictly positive  $\delta_0$ exists, then 
the  probability $\P \left\{ \Vert A^{(n)} \Vert > 2v(1+xn^{-2/3}) \right\}, x>0$ is 
bounded in the limit of infinite $n$ by an expression that does not depend on the
details of the probability distribution of $a_{ij}$.  
The proof is based on the completed and modified version of the approach 
developed by \mbox{Ya. Sinai} and A. Soshnikov to study high moments of Wigner random matrices.

\end{abstract}

\section{Introduction}

Random matrices of infinitely increasing dimensions are of 
the intensive study 
due to their deep mathematical content and numerous applications 
mathematical and theoretical physics and in various branches of applied mathematics
(see review  \cite{B} and monograph \cite{M} and references therein). This  can be also explained by fairly 
wide  universality of spectral properties of random matrix ensembles,
in particular, the asymptotic universality of the local eigenvalue statistics 
in the bulk and at the border of the limiting spectra predicted first in the early physics studies. 

\vss
The spectral theory of large random matrices has been started by the work of \mbox{E. Wigner} \cite{W}, 
where the limiting eigenvalue distribution
of   the ensemble of $n\times n$ random real symmetric matrices of the form
$
A^{(n)}_{ij} =  a_{ij}/\sqrt n
$
has been  considered. 
Random variables $\{a_{ij}\}_{1\le i\le j\le n}$ are jointly independent and have symmetric probability distribution.
Assuming that  the variance of all  random variables $a_{ij}$ is equal to $v^2$ and 
that all  moments of $a_{ij}$ exist, E. Wigner has proved 
the convergence 
$$\lim_{n\to\infty} \  {1\over n} \E L_{l}^{(n)} =  m_{l} = 
\cases{
v^{2k} {(2k)!\over k! (k+1)!},  & if $ \   l=2k$, \cr
0, & if \ $ l=2k+1$,\cr
}
\eqno (1.1)
$$
where
$$
L_{l}^{(n)} = 
\T (A^{(n)})^{l} = {1\over n} 
\sum_{j=1}^n \left(\lambda_j^{(n)}\right)^{l}, \quad l\ge 0 
\eqno (1.2)
$$ 
with real eigenvalues  $\l_j^{(n)} = \l_j(A^{(n)})$, $\l_1^{(n)}\le \dots\le \l_n^{(n)}$. Here and below $\E$ denotes the 
mathematical expectation with respect to the family $\{a_{ij}\}$.  \mbox{E. Wigner}  
has shown that the family $\{m_{l}, l \ge 0\}$ represents the moments of the measure with the density of the  form
$\sqrt{ 4v^2- \l^2}/({2\pi v^2})$
 and the support $[-2v,2v]$ \cite{W}. 
This measure is referred to as to the semicircle distribution and  convergence (1.1) 
reflect the measure convergence  known as the Wigner or the semicircle law \cite{G,P}. This convergence reflects the 
properties of the eigenvalue distribution of $A^{(n)}$ that can be referred to as to the global ones.

\vss

In the particular case, when $a_{ij}$ are given by real jointly independent \mbox{Gaussian}
random variables,
the family $\{A^{(n)}\}$  is referred to as to the Gaussian Orthogonal Ensemble abbreviated as GOE; 
the corresponding ensemble of hermitian random matrices is known as the Gaussian Unitary Ensemble (GUE) \cite{M}.
In these cases the variance of the diagonal elements $a_{ii}$ 
is slightly  different  with respect to those of the  random matrices studied by E. Wigner.
 
The famous Tracy-Widom law says that the probability distribution of the maximal eigenvalue of GOE matrices converges,
when rescaled, to a probability distribution known as the Tracy-Widom distribution \cite{TW}
$$
\lim_{n\to\infty} \P\left\{ \l_{\max} ^{(n)} > 2v \left(1+  {x\over n^{2/3}}\right) \right\}  = F_{TW}^{(1)}(x),
\eqno (1.3)
$$
where the spectral norm  $ \l_{\max}^{(n)} = \max\{ \vert \l_1^{(n)}\vert ,\vert \l_n^{(n)}\vert \}=  \Vert A^{(n)}\Vert$.
The same result is valid for the maximal eigenvalue  of GUE matrices with $F_{TW}^{(1)}(x)$ 
replaced by another function $F_{TW}^{(2)}(x)$. The statements like (1.3) 
that concern the finite number of eigenvalues of $A^{(n)}$ can be referred to as to the local ones.

There exist a large number of results showing  strong and deep connections between 
the local eigenvalue distributions of the form (1.3)  of random matrices from one side and 
the probability distribution of a number of combinatorial and probabilistic objects 
(see for example papers and reviews \cite{Baik,Boug,J,Ok} and references therein). Thus
the universality of the local spectral properties   can be considered as one of the most important 
 issue in the spectral theory of random matrices.

 \vss The local properties of large random matrix spectra can be studied with the help of (1.2) in the limit
 when $l$ and $n$ simultaneously tend to infinity  \cite{BY,FK,G}. 
In papers \cite{SS1,SS2}, a  powerful approach has been proposed and elaborated  
by Ya. Sinai and A. Soshnikov
to study the mean values of the traces $ L_{2s}^{(n)}$ of  $A^{(n)}$    and their 
correlation functions in the 
limit $s,n\to\infty$.
This approach is based on the representation of the mean 
value $M_{2s}^{(n)}= \E L_{2s}^{(n)}$
as the sum over weighted trajectories $I_{2s}$  of $2s$ steps that can be classified according to the number of their 
self-intersections.
The upper bounds for $M_{2s}^{(n)}$ and the central limit theorems for $L_{2s}^{(n)}$ were established 
in these papers
for the asymptotic regimes when  $s= o(n^{2/3}), n\to\infty$. 

In the subsequent  paper \cite{S}, the  asymptotic regime   $s= O(n^{2/3})$ necessary to prove the statements of the form (1.3)
has been reached.  In \cite{S}, the limiting correlation functions of $L_{2s}^{(n)}$ 
have been studied 
under hypothesis that  the moments of $a_{ij}$ are bounded as follows
$$
V_{2m} =  \sup_{1\le i\le j} \E (a_{ij})^{2m} \le (C_a \, m)^m \quad {\hbox{for all \ }} m \in {\bf N}.
 \eqno (1.4)
 $$ 
To relax the  condition (1.4), A. Ruzmaikina \cite{R} has further
specified  the method proposed by Sinai and Soshnikov.
In \cite{R}, the polynomial decay of the probability distribution of $a_{ij}$ is assumed such that
$$
\sup_{1\le i\le j} \P\{ \vert a_{ij}\vert >y \} \le C_b\, y^{-18} .
\eqno (1.5)
$$ 
It was also indicated  that  in paper \cite{S}, a certain part of the sum  that correspond to  trajectories 
$I_{2s}$ with large number of steps out from the same site 
has not been correctly estimated and an more precise computations  has been presented in \cite{R}.  
\vss

To aim of the present paper is to study the moments $M_{2s}^{(n)}$ in the limit $n\to\infty$, $s = O(n^{2/3})$
under  more weaker   conditions than that of (1.5). To do this,  
we introduce 
a new modified version of the approach by Ya. Sinai and A. Soshnikov improved by A. Ruzmaikina.  
However, this approach in its original version is not complete \cite{KV}. 
The doubtful point of the work \cite{R}
concerns again the properties of the   trajectories
mentioned above that have many steps out from the same site. 
Following \cite{S}, it is claimed  that 
the presence of such sites  can happen in two situations. The first one arises 
 when  the underlying Dyck path
(or equivalently,  the related plane rooted tree) contains vertices of high degree; 
in the opposite case, these sites can be created by
 the arrivals to them along the ascending steps of the Dyck paths.

The third possibility, complementary to the two first ones,
has not been detected in \cite{R,S} and the analysis 
does not cover the whole set of trajectories.
This  third   situation is given by the case when  the trajectory  creates 
the sites with large number of exit steps 
after arrivals at these sites by  some of descending steps of the corresponding Dyck path.
This special kind of arrivals can happen  when the trajectory has previously performed a  number of 
fractures of the tree structure of the walk. 
In paper \cite{KV}, we have studied in details
this type of walks with broken tree structure and have  completed the Sinai-Soshnikov technique. 
The example case of random matrices with bounded random variables has been considered
in \cite{KV}; also a way to study  the case when $a_{ij}$ has a finite number of 
moments has been indicated there.

\section{Main results and scheme of the proof}

Let us consider the   ensemble $\{A^{(n)}\}$ of $n\times n$ real symmetric random matrices with  elements
$$
A^{(n)} = {1\over \sqrt n} a_{ij}, \quad 1\le i\le j\le n,
\eqno (2.1)
$$
where the family of jointly independent identically distributed random \mbox{variables} 
\mbox{${\cal A} =  \left\{ a_{ij}, 1\le i\le j  \right\}$} 
determined on the same probability space $(\Omega, \P, {\cal F} )$ is such that 
$$
\E a_{ij} = 0 \quad {\hbox{and \ \ }} \E a_{ij}^2 = v^2.
\eqno (2.2)
$$
Here and below we denote by $\E$ the mathematical expectation with respect to the measure $\P$. 
We refer to the ensemble (2.1), (2.2) as to the Wigner ensemble of random matrices.
We denote $
V_{2m} =  \E  a_{ij}^{2m}$, 
$2\le m\le 6 $. Our main result is as follows.

\vs 
\noindent {\bf Theorem 2.1.} {\it Consider the Wigner random matrix ensemble (2.1), (2.2) and assume that the random variables 
$a^{(n)}_{ij}$ have symmetric probability distribution and that there exists $\delta_0>0$ such that 
$$
V_{12+2\d_0} =  \E \vert a_{ij}\vert^{12+2\d_0} < +\infty.
\eqno (2.3)
$$
Then for any $x>0$, 
$$
\limsup_{n\to\infty} \P \left\{ \l_{\max}(A^{(n)}) > 2v\left(1+{x\over n^{2/3}}\right)\right\} \le \inf_{\chi >0} \CL(\chi) e^{- x\chi },
\eqno(2.4)
$$
where $\CL(\chi)$ does not depend on the 
particular values of $V_{2m}, 2\le m\le 6$; this value is determined as the following limit of the moments of GOE}
$$
\CL (\chi) = \lim_{n\to\infty}\  {1\over (2v)^{2\lfloor s_n\rfloor}} \E \T ( A^{(n)}_{\hbox{\tiny{GOE}}})^{2\lfloor s_n\rfloor},
\quad s_n = \chi n^{2/3},
\eqno (2.5)
$$
{\it 
where $\lfloor s\rfloor$ is the largest integer not greater than $s$.
}

\vs
Existence of the limit (2.5) is established in \cite{S}.
Theorem 2.1 is a consequence of the following statement that represents the main technical result of the paper.

\vs 
\noindent {\bf Theorem 2.2.} {\it 
Given $0<\d<1/6$, consider the truncated random variables
$$
\hat a_{ij} ^{(n)} = \cases{
 a_{ij},  & if$\ \vert a_{ij} \vert \le U_n = n^{1/6 - \delta}$;\cr
0, & otherwise,
\cr
}
\eqno (2.6)
$$
such that  $a_{ij}$ verify conditions of Theorem 2.1. 
Then  the high moments of the random matrices $\hat A^{(n)}$ with the elements 
$\hat a^{(n)}_{ij}/\sqrt n$ verify relation
} 
$$ {\cal L}(\chi) = \lim_{n\to\infty} {1\over (2v)^{2\lfloor s_n \rfloor} }  \E  \T ( \hat A^{(n)})^{2\lfloor s_n\rfloor} 
<+\infty,
\eqno (2.7)
$$
{\it
where  $s_n = \chi n^{2/3}$ and $\chi>0$. 
The value of $ {\cal L}(\chi)$ (2.7) does not depend on particular values of the moments $V_{2m}, \, 2\le m\le 6$
and therefore coincides with that of (2.5).
}

\vs
{\it Remarks. } 

1. The main inequality we prove is of the following form that coincides 
with the expressions obtained previously \cite{S},
$$
\limsup_{n\to\infty} {1\over v^{2\lfloor s_n \rfloor} }  \E \T ( \hat A^{(n)})^{2\lfloor s_n\rfloor} 
\le {1\over \sqrt{\pi \chi^3}}\,  B(6\chi^{3/2})\, e^{C\chi^3},
\eqno (2.8) 
$$
 where $B(\tau)$, $\tau>0$ is determined 
 by relation 
$$
B(\t) = \lim_{k\to\infty}  B_k(\tau), \quad B_k(\t) = {1\over \vert \Theta_{2k}\vert} \ \sum_{\theta\in \Theta_{2k}} 
\exp\left\{
{\t\over \sqrt {k} }\  \max_{ 1\le l \le 2k}  \theta(l) \right\},
\eqno (2.9)
$$
where $\Theta_{2k}$ is the set of all $2k$-step Dyck paths, i.e. the 
simple walks $\theta$ of $2k$ steps  that start and end at zero 
and stay non-negative. It is known   \cite{St}
 that the cardinality  $\vert \Theta_{2k}\vert$ is given by the Catalan number
that we denote by $t_k= (2k)!/k! (k+1)!$. In what follows, we will use also
the variable $\tilde B_k(\tau)  = \inf_{k\ge k'} B_{k'}(\t)$. 

The bound of the form (2.8)  
has been derived  for the first time  in 
\cite{SS2};  existence of the limit  (2.9) is rigorously proved in  \cite{KM}. 
It is shown  in \cite{KM} that   $B(\t)$  coincides with the corresponding exponential moment
of the normalized Brownian excursion on $[0,1]$.



 \vss
  
 2. Theorems 2.1 and 2.2   remain true in the case when the Wigner ensemble of hermitian random matrices is considered
 instead of the real symmetric ones.

\vss Let us  describe  the scheme of the proof of Theorem  2.2. 
The general \mbox{approach} of \cite{SS1,SS2} and \cite{R,S}  is based on the Wigner's original  point of view 
when the trace of the power of  $\hat A^{(n)}$ 
is  represented as the weighted sum 
$$
\E  \T (\hat A^{(n)})^{2s} = {1\over n^s} \sum_{i_0, i_1, \dots, i_{2s-1} =1}^n \E \left\{ \hat a_{i_0i_1}
\cdots \hat a_{i_{2s-1}i_0}\right\} = {1\over n^s} \sum_{I_{2s} \in \CI_{2s}(n)} \hat \Pi (I_{2s}), 
\eqno (2.10)
$$
where the sequence $I_{2s}=  (i_0, i_1,  \dots, i_{2s-1},i_0)$,  $ i_l \in \{1, 2,\dots, n\} $ is regarded as a closed trajectory of $2s$ steps,
 and the weight $\Pi(I_{2s})$ is given by the average value of the product of corresponding random variables
 $a_{ij}$. It is natural to say that the subscripts of the variables $i_t$ represent the values of the discrete time from the interval 
 $[0,\dots, 2s]$. 

Assuming that the probability distribution of $a_{ij}$ is symmetric, 
the non-zero contribution to the right-had side of (2.10) comes from the trajectories $ I_{2s}$
whose weight $\hat \Pi( I_{2s})$ contains each random variable $a_{ij}$ with  even multiplicity. 
In \cite{SS1}, these trajectories were referred to  as to  the even
ones. Regarding the even walks,  
one can define the marked instants of time as the instants,  when each
random variable is seen for the odd number of times. 

The principal idea of the Sinai-Soshnikov
method is to separate the set of all possible even trajectories $ \CI_{2s}$ into the classes of equivalence
according to the numbers of vertices of the self-intersections of given degree $\kappa$. 
A. Ruzmaikina has shown that the vertices  of self-intersections of high degree 
$\kappa>k_0$ can be neglected, where $k_0$ depends on the order of the maximal existing  moment
of the random variables $a_{ij}$. 

\vss

In the present paper, we further modify
the technique developed 
by Ya. Sinai and A. Soshnikov by introducing  more informative classes of equivalence.
In some part, our  representation is inspired by the ideas of the earlier work of  Z. F\"uredi and J. Koml\'os \cite{FK}.

\section{Walks, graphs and classes of equivalence}

In the present section we introduce classes of equivalences of trajectories $I_{2s}$ (2.10)  and 
develop a method to estimate the cardinalities of these classes.
In subsection 3.1, we give the definitions of walks and graphs. 
In subsection 3.2, we describe the Sinai-Soshnikov classes of equivalence.
Further modification of  the Sinai-Soshnikov classes
is presented in subsection 3.3.
In subsection 3.4,  the number of walks in these new classes of equivalence is estimated.

\subsection{Trajectories, walks and graphs of walks}

Regarding  $I_{2s}$ (2.11), we write that $I_{2s}(t) = i_t$ for integer  $t\in [0,2s]$ and say that the couple
$(i_{t-1},i_t)$ represents the step number $t$ of the trajectory $I_{2s}$. Let us consider 
the set 
$
\CU(I_{2s};t) = \{ I_{2s}(t'), \ 0\le t'\le t\}
$
and denote by $\vert \CU(I_{2s};t)\vert $ its cardinality.
Given a particular trajectory $I_{2s}$, we construct corresponding  {\it walk } 
$w(I_{2s};t)= w_{2s}(t)$, $0\le t\le 2s$  by the \mbox{following} recurrence rules: assuming  that there is
an infinite ordered alphabet $\{\a_1, \a_2, \dots\}$, we say that 
\vsss
\noindent 1) at the initial instant of time, $w_{2s}(0) = \a_1$;

\noindent 2) if $I_{2s}(t+1) \notin \CU(I_{2s};t)$, then $w_{2s}(t+1) = \a_{\vert \CU(I_{2s};t)\vert +1}$;

 if there exists such $t'\le t$ that $I_{2s}(t+1) = I_{2s}(t')$, then $w_{2s} (t+1) = w_{2s}(t')$.

\vsss \noindent One can say that $I_{2s}$ generates $w_{2s} = w(I_{2s})$.

Certainly, there exists a number of different  trajectories $I_{2s}$ that generate the same  walk. 
For example,  $I'_6 = (5,2,3, 5,2,3,5)$ and $I''_6= (3,7,8,3,7,8,3)$
are such that  $w(I'_6) = w(I''_6) = (\a_1,\a_2,\a_3,\a_1,\a_2,\a_3,\a_1)$. 
Then one can determine 
a partition of  the set $\CI_{2s}$ of all possible trajectories  of $2s$ steps  into classes of equivalence 
labelled by the walks  $w_{2s}$. We denote the corresponding classes of equivalence by  
$\CI_{w_{2s}}$. 
Regarding the walks $w_{2s}$, we will  use sometimes,  instead of the  alphabet $\{\a_1, \a_2,\dots \}$, 
 the greek letters  without subscripts.
 The first symbol of the walk is called the {\it root} of the walk.


\vss

Given  $w_{2s}$, it is useful to consider a natural  {\it graphical representation} of the walk  $\Gamma(w_{2s})= (\CV, \CE)$, 
where $\CV = \CV(w_{2s})$
is a set of vertices labelled by the symbols (or  letters)  from the alphabet. 
The set of oriented ordered edges $\CE(w_{2s})$ 
contains a couple $e=(\a,\b)$, $\a,\b\in \CV(w_{2s})$, if there exists an instant $t$ such that $w_{2s}(t-1)=\a$ and $w_{2s}(t)=\b$. 
In this case we write that $(\a,\b) = e(t)$.
Clearly,  $\vert \CE(w_{2s})\vert = 2s$. 
The root of the walk determines the {\it root vertex}
of $\Gamma$. We denote it by $\rho$  in the context where its difference with other vertices is important.

Given a walk  $w_{2s} $, one can also determine the 
graph $ \CF(w_{2s})= (\CV,\tilde \CE)$
with non-oriented simple edges  $ e_i \in \tilde \CE$, $1\le i\le s$ that we denote by  $\{\a,\b\}$.  
We say that $\{\a,\b\} \in \tilde \CE$ if there exists at least one edge
of the form $(\a,\b)\in \CE$ or $(\b,\a)\in \CE$.  We refer to the graph  $\CF(w_{2s})$ as to the {\it frame of the walk}
$w_{2s}$.
\vss

We  determine the 
{\it current multiplicity} of  the frame edge $\{\a,\b\}$ at  \mbox{time $t$}
by variable
\begin{eqnarray*}
\CM_{w}(\{\a,\b\};t)  = \# \{ t'\in [1,t]: \  (w(t'-1), w(t')) &=& (\a,\b) \ {\hbox{or}}\  \\ 
 (w(t'-1), w(t')) &=& (\b,\a)\}.
\end{eqnarray*}
The probability law of random variables $a_{ij}^{(n)}$ being symmetric, all their odd moments
vanish. Then the non-zero contribution to (2.10) comes from the closed trajectories such that 
for each frame edge $\{\a,\b\}$ equality 
$
\CM_{w}(\{\a,\b\};2s) = 0({\hbox{mod }}2)
$ is true.
We refer to 
such walks as to the {\it even  walks}.  We denote by $\CW_{2s}$ 
the set of all possible even  walks of $2s$ steps. In what follows, we consider the even  walks only
and refer to them simply as to the walks.
\vsss

Given $w_{2s} \in \CW_{2s}$, we say that the instant of time $t$ with $w_{2s}(t) = \b$ is {\it marked} 
if the frame edge $\{\a,\b\}$ with $\a = w_{2s}(t-1)$ is passed an odd number of times during the interval $[0,t]$;
$
\CM_w(\{\a,\b\}; t) = 1({\hbox{mod }} 2).
$
In this case we will say that the step $t$ of the walk and the corresponding edge $e(t)\in \CE$ of $\Gamma(w_{2s})$ are marked. 
All other instants of time, steps and edges are referred to as the {\it non-marked} ones. 
Clearly, the number of the marked steps of the even walk is equal to the number of non-marked steps.

Regarding the graphical representation $\Gamma(w_{2s})$, we can  remove all non-marked edges from the set $\CE(w_{2s})$
and keep the set of all marked edges $ \CE'(w_{2s})$. Then we get a new graphical representation 
$\Gamma'(w_{2s}) = (\CV,\CE')$, where the oriented edges  $e'\in \CE'$ are ordered according their time labels due to $w_{2s}$.

We say that $M= \CM_{w}(\{\a,\b\};2s)$ is the number of times that $w_{2s}$ passes the frame edge 
$\{\a,\b\}\in \CF$.   Also we say  that the number $M/2$ determines
the {\it multiplicity} of the edge $\{\a,\b\}$. In the case when 
$M/2>1$ we will say that 
$\{\a,\b\}$ is {\it multiple} or {\it $M/2$-fold edge}. It is clear that the multiplicity of the frame edge 
$\{\a,\b\}$ coincides with the number of elements $e'\in \CE'$ such that $e' = (\a,\b)$ or $e'=(\b,\a)$. 
Therefore, with some language abuse,  
we can simply say that $M/2$ represents the multiplicity of the edge $(\a,\b)$ of the walk $w_{2s}$.


\begin{figure}[htbp]
\centerline{\includegraphics[width=12cm,height=4cm]{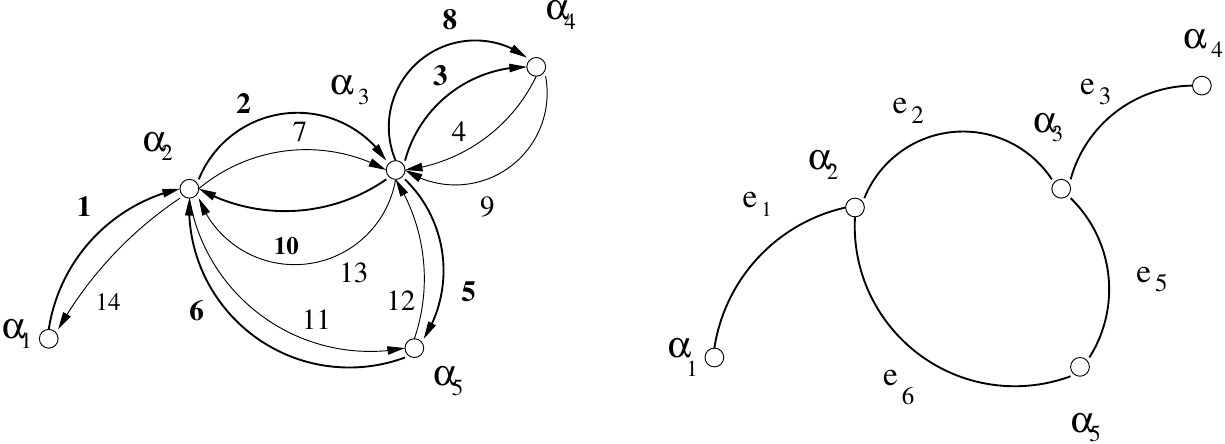}}
\caption{\footnotesize{ The graph   $\Gamma(w_{14})$ of the walk $w_{14}$ and its frame
$ \tilde \Gamma$}}
\end{figure}

\vsss
Given $\b\in \CV$, we determine its {\it exit cluster}  ${\Delta}(\b)$ as the set of marked edges $e'(t_i)\in \CE'$
such that $w_{2s}(t_i-1)=\b$. Sometimes we will use $\Delta(\b)$ to denote the set of corresponding vertices.
We determine the {\it exit degree} of $\b$ as 
$ \vert {\Delta}(\b)\vert$ and denote 
the maximal exit degree of the walk $w_{2s}$ as
$$
D(w_{2s}) = \max _{\b \in \CV}  \vert {\Delta}(\b)\vert.
$$
In what follows, 
we will refer to  the graphical representation $\Gamma'(w_{2s})$ 
simply as to the multigraph, or simply the {\it graph  of the walk}  $w_{2s}$.
Let us stress that  that all elements of $\CE'$ are ordered, or in other words, have their  time labels.

\vsss
On Figure 1,  we have shown the graphical representation  $\Gamma_{14}= \Gamma(w_{14})$ of the walk 
$$
w_{14} = (\a_1,\a_2,\a_3,\a_4,\a_3,\a_5,\a_2,\a_3,\a_4,\a_3,\a_2,\a_5,\a_3,\a_2,\a_1)
$$
and as well as its frame  $\CF(w_{2s})$. The marked edges of $\Gamma_{14}$ are given in boldface. 
The corresponding graph $\Gamma' = (\CV,\CE')$ is such that its edges 
$e'\in \CE'$ have time labels $\{1,2,3,5,7,8,10,11\}$. There are two double edges in the walk $w_{14}$, these are 
$(\a_2,\a_3)$ and $(\a_3,\a_4)$. The maximal exit degree of $w_{14}$ is $D(w_{14})=3$.

\subsection{Instants of self-intersection and Sinai-Soshnikov classes }

Any closed even walk $w_{2s}\in \CW_{2s}$ generates a binary sequence 
$\theta_{2s} = \theta(w_{2s})$ of $2s$ elements $0$ and $1$ that corresponds to the non-marked and the marked steps of $w_{2s}$s,
respectively. The set of such sequences is in one-by-one correspondence with the set of the Dyck paths of $2s$ steps 
\cite{St}. We denote by $\Theta_{2s}$ the set of all Dyck paths of $2s$ steps and say that 
$\theta_{2s} = \theta(w_{2s})$ is the Dyck structure of $w_{2s}$.

Regarding $\b\in \CV(g_{2s})$, we denote by $1\le t_1^{(\b)}\le \dots\le t^{(\b)}_N$ all marked instants of time
such that $w_{2s}(t_j^{(\b)}) = \b$. We say that the $N$-plet $(t^{(\b)}_1,\dots, t^{(\b)}_N)  $ 
represents the {\it arrival marked instants} at $\b$
and that $e(t_j^{(\b)})$ is the {\it arrival marked  edge}  at $\b$. 
Following \cite{SS1}, we say that $\b$ is the vertex of $N$-fold self-intersection.

For any non-root vertex $\b\in \CV(g_{2s})$,  $\b\neq \rho$, we have $N= N_\beta \ge 1$. 
We define also  the variable $N^{(t)}_\b$ given by the number of 
marked arrival instants at $\b$ during the time interval $[0,t]$.

\vsss
If $N_\b=2, \b\neq \rho$, then the vertex $\b$ is said to be the {\it vertex of  two-fold (or simple) self-intersection}
and $t^{(\b)}_2$ represents the {\it instant of simple self-intersection} \cite{SS1}. We will also say that the couple
$(t_1^{(\b)},t_2^{(\b)})$ represents the two-fold self-intersection of the walk.
\vsss
If $N_\b = k, \b\neq \rho$, 
then $\b$ is said to be the vertex of $k$-fold self-intersection and $k$ is the self-intersection degree of $\b$.
In general, we denote by $\kappa = \kappa(\b)$ the {\it self-intersection degree} of $\b$.
For the root vertex $\rho$, we  set $\kappa(\rho)=N_\rho+1$.

\vss
The following definition is very important in the studies of the even closed walks.
\vsss
\noindent {\bf Definition 3.1. \cite{SS2}} {\it The arrival  instant $t=t^{(\b)}$ with $w_{2s}(t) = \b$  
is said to be the non-closed (or open) instant of simple self-intersection
if the step $(t-1,t)$ is marked and if there exists
at least one non-oriented edge $\{\b,\g\}\in \tilde \CE$ attached to $\b$ that
is passed an odd number of times during the time interval $[0,t-1]$. In this case we say that 
the edge $\{b,\g\}$ is open up to the arrival instant $t=t^{(\b)}$, or more briefly that this edge is $t$-open. 
Also we can say in this case  that the vertex $\b$ is $t$-open vertex of simple self-intersection.
}
\vsss
{\it Remark.}  
Definition 3.1 remains valid in the case when $\g=\b$ and the graph $g_{2s}$ has a loop at the vertex $\b$, i.e. the edge $(\b,\b)$. 
For example, the walk $(\a,\b,\b,\a)$ has a simple open self-intersection $(1,2)$.  We assume that the walk
$(\a,\a,\a)$ has a simple self-intersection that is not open. 

One more example is given by the walk 
\mbox{$w_6 = (\a,\b,\g,\a,\b,\g,\a)$,} where the root vertex  is the vertex of the two-fold open self-intersection, $\kappa(\a)=2$.

\vss
We say that the walk $w_{2s}\in \CW_{2s}$ belongs to the class  $\CC^{(\th)}{(\bar \nu,r)}$ \cite{SS2},
where the set $\bar \nu = \bar \nu(w_{2s}) = ( \nu_2, \nu_3,\dots,\nu_s)$, 
when it has the sequence of marked edges determined by $\th$, 
and the graph of the walk contains  $\nu_2$ vertices of two-fold self-intersections,
$\nu_3$ vertices of three-fold self-intersection, etc.; among them there are $r$ open vertices  of two-fold self-intersection.
In papers \cite{SS2,S} the classes of walks that have multiple edges have been studied
and the upper bounds for the moments of random matrices are obtained under conditions (1.4). To relax these conditions,
one has to study  the walks with multiple edges in more details. 
This can be done 
by introducing a classification more informative than that given by $\CC^{(\theta)}(\bar \nu,r)$.

\subsection{Classes of equivalence $\CC_{k_0}(\bar \mu,P, \bar Q; \bar \nu)$}

Let $k_0$ be a given integer. 
Looking at the graph of the walk $\Gamma'(w_{2s})$ that contains the marked edges only, 
we consider the vertices $\a_i$ such that $\kappa(\a_i)\ge k_0+1$
and say that they are the \mbox{ {\it $\nu$-vertices}.} The remaining ones are referred to as to the {\it $\mu$-vertices}. 
Regarding a $\nu$-vertex $\a$, we say that all edges of the form $(\b_j,\a) \in \CE'$   are the {\it $\nu$-edges}.

If $\a$ is a $\mu$-vertex joined with another $\mu$-vertex $\b$ by one or more edges, 
one can consider the minimal edge of the form $\{\a,\b\}$, for instance $(\a,\b)$, and say that 
$(\b,\a)$ is the $\mu$-edge. Under the minimal edge we mean here the edge whose time label is less than 
that of any other one of the form $\{a,\b\}$. We say that the next to the minimal edge of the form
$\{a,\b\}$ is the {\it $p$-edge}. The next after the  $p$-edge element $\{\a,\b\}$ is said to be the {\it $q_1$-edge}.
The {\it $q_j$-edges} are determined by recurrence. 

\vss
If the $\mu$-vertex $\a$ is joined with a $\nu$-vertex $\b$, then the minimal edge of the form
$(\b,\a)$, if it exists,  is referred as to the $\mu$-edge. Then one can determine the $p$-edge $(\b,\a)$, and so on. 
One can think about the $\nu$-edges as those colored in black, the $\mu$-edges colored in blue
and the $p$-edges colored in green. The $q$-edges remain non-colored. 
Regarding two $\mu$-vertices $\a$ and $\b$ that are joined by edges, we say that the $p$- and $q$-edges
of the form $\{\a,\b\}$ are attributed to the blue $\mu$-edge that joins $\a$ and $\b$. 

\vss
Regarding a $\mu$-vertex $\a$ that is not the root vertex $\rho$, we determine the degree of $\mu$-self-intersection of $\a$ 
as the total number $m$ of $\mu$-edges of $\Gamma'$ that arrive at $\a$ and denote it by $\kappa_\mu(\a)=m$. 
As for the root vertex, we denote by $m_\rho$ the number of $\mu$-edges of the form $(\a_i,\rho)$
and set $\kappa_\mu(\rho) = m_\rho+1$. Certainly, all these definitions remain valid
in the case when there are the loop $\mu$-edges, as well as the loop $\nu$-edges.

\vss

A  walk $w_{2s}$ belongs to the 
 class  $\CC_{k_0}(\bar \mu,  P, \bar Q;\bar \nu)$,
with {$\bar \mu = (\mu_1,\mu_2,\dots, \mu_{k_0})$,} $\bar Q = (Q_1, Q_2,\dots, Q_{k_0-2})$
and $\bar \nu = \bar \nu^{(k_0)} = (\nu_{k_0+1}, \dots, \nu_s)$, if its graph $\Gamma'$
contains $\mu_m$ vertices with the $\mu$-self-intersection degree $m$, $P$ $p$-edges, 
$Q_j$ $q_j$-edges and $\nu_k$ vertices of the self-intersection degree $k\ge k_0+1$. 

\vss

Let us consider Figure 1 and assume that $k_0\ge 3$. Then all vertices of $\Gamma'$ are the $\mu$-vertices
and therefore the walk $w_{14}$ belongs to the class $\CC_{k_0}(\bar \mu, P,\bar Q;\bar \nu)$
with $\mu_1= 4$ (including the root vertex), $\mu_2=1$, $P=2$, $\bar Q = \bar 0$ and $\bar \nu = \bar 0$.


\subsection{Estimates of the number of walks}

Any walk $w_{2s}$ is determined by its position (or the alphabet value) 
at the instants of time $t$, 
$0\le t\le 2s$.  Regarding the marked instants of time, we see that $w_{2s}$ generates a partition of them into groups,
each group corresponding to the vertices of $\Gamma(w_{2s})$ that are determined by the
first arrival instants at them. Then one can get the upper bound for the number of walks in the same class
of equivalence by estimating the number of possibilities to split the set of marked instants of time into groups.
This is the principal idea of the method developed by Ya. Sinai and A. Soshnikov. 
 Let us describe this technique  in more details that we will also need in our reasoning.

\subsubsection{Sinai-Soshnikov estimates}

Any Dyck path $\th \in \Theta_{2s}$ generates an ordered sequence of $s$ marked instants of time 
$1= \xi_1<\xi_2<\dots <\xi_s\le 2s-1$ that we 
 denote by $\Xi_s = (\xi_1, \dots,\xi_s)$.

\vss 

Let us denote by $\Psi_s(\bar \nu)$  the  number of possibilities to produce a partition of the set $\Xi_s$ of $s$ labelled 
elements that contains  $\nu_k $ non-intersecting subsets of $k$ elements, $k=1,\dots, s$
such that $s= \sum_{k=1}^s k\nu_k$. Then it is not hard to see that
$$
\Psi_s(\bar \nu) =  {s!\over (s-\sum_{k\ge 2} k\, \nu_k)! \ (2!)^{\nu_2}\,  \nu_2! \ (3!)^{\nu_3}\,  \nu_3!\  \cdots\  (s!)^{\nu_s}\,  \nu_s!}
\le \prod_{k=2}^s {1\over \nu_k!}  \left( {s^k\over k!}\right)^{\nu_k},
\eqno (3.1) 
$$
where we have used  inequality ${s!/ (s-a)!} \le s^a$ with $a= \sum_{k\ge 2} k\nu_k$.
One can \mbox{consider} the last product over $k$ (3.1) as the product over the vertices of self-intersection.
This gives a prescription for the walk 
where to go at the next-in-turn marked instant of time.
No such a prescription  is given for the non-marked instants of time. 
If the walk leaves the vertex of self-intersection $\b$ at the non-marked instant, 
it can have a number of vertices where to go 
 bounded  by certain value $\Upsilon(\b)$.

To get the upper bound  for the number of all possible walks,  we have to multiply the right-hand side of (3.1)
 by the upper bound $\Upsilon$. It is shown in \cite{SS1} that if $\kappa_\nu(\b) = k$, $k\ge 2$, then 
 $\Upsilon(\b) \le  \Upsilon_k =  (2k)^k$. For the moment, we ignore in (3.1) the 
 special role of the root vertex (see Appendix C, subsection 5.3  for more details).  Multiplying  the 
 right-hand side of (3.1) by $\prod_{\b} \Upsilon(\b)$, we get the upper bound for the number of walks
 in the class of equivalence $\CC(\bar \nu)$.

 Given $\theta\in \Theta_{2s}$, let us consider the 
 class  $\CC^{(\theta)}_d(\bar \nu,r,p)$ of corresponding to $\theta$  walks $w_{2s}$
 such that the maximal exit degree in $\Gamma'(w_{2s})$ is equal to $d$, 
 $r$ is the number of open simple self-intersections and $p$ is the number of simple self-intersections
that produce the  edges of multiplicity two with the same orientation of the marked edges.
In
 \cite{SS2,S}, it is shown  that the number of walks in this class 
 is bounded as follows
$$
\vert \CC^{(\theta)}_d(\bar \nu,r,p)\vert  \le {1\over (\nu_2-r-p)!} \left( {s^2\over 2} \right)^{\nu_2-r-p}
\cdot {(6sH_\theta)^r\over r!} \cdot {(3sd)^p\over p!} 
$$
$$
\times  \prod_{k=3}^s \,  {1\over \nu_k!} \left( {s^k  \, \Upsilon_k\over k!} \right)^{\nu_k} ,
\eqno (3.2)
$$
where $H_\theta = \max_{1\le t\le 2s} \theta(t)$ and  $ \U_k= (2k)^k $.  When deriving (3.2) in \cite{SS2,S},
some  arguments additional  to the combinatorial reasoning of (3.1) are used.
In Section 5,  we describe the proof  (3.2)  in slightly more general setting.

\subsubsection{Estimate for the classes $\CC(\bar \mu,P,\bar Q;\bar \nu)$}

\vss
Given $\theta\in \Theta_{2s}$ and $k_0\ge 2$, let us consider 
the class  of the walks
$\CC^{(\th)}_{d,k_0}(\bar \mu,  P, \bar Q,  r;\bar \nu )$ such that the maximal exit degree of the corresponding
$\Gamma'$ is equal to $d$ and  $r$ is the number of  
simple self-intersections  that are the open ones.

  \vs

{\bf Lemma 3.1.} {\it
 If 
$$
\sum_{m=2}^s (\mu-1)\mu_m + P + \vert \bar Q\vert + \sum_{k=k_0+1} (k-1)\nu_k \le {s-1\over 6},
\eqno (3.3)
$$
where
 $\vert \bar  Q\vert = \sum_{j=1}^{k_0-2} Q_j$, 
 then  the cardinality of  the class $\CC^{(\th)}_{d,k_0}(\bar \mu,  P, \bar Q,  r;\bar \nu )$ is bounded as follows (cf. (3.2)),
  $$
  \vert \CC^{(\th)}_{d,k_0}(\bar \mu,   P, \bar Q, r;\bar \nu )\vert 
  \le  {1\over (\mu_2-r)!} \left( {s^2\over 2} \right)^{\mu_2-r}\cdot 
  {(2sH_{\th})^r\over r!} 
  $$
  $$
  \times   \prod_{m=3}^{k_0}\,  {1\over \mu_m!} \left( { s^m\over m!} \right)^{\mu_m}
   \cdot \,  { (2s \tilde d )^{P}\over P!} 
  \cdot   {(2\tilde d P)^{Q_1}\over Q_1!}   \ \ 
  \prod_{j=2}^{k_0-2}  {(2\tilde d Q_{j-1} )^{Q_j}\over Q_j!} 
  \ \ \prod_{k=k_0+1}^s {1\over \nu_k!}   \left( { s^k\over k!}\right)^{\nu_k}
  $$
  $$
  \times  \U_{k_0}(\bar\mu, P, \vert \bar Q\vert, r;\bar \nu ),
  \eqno (3.4)
  $$
with $\tilde d = \max\{ d,k_0\}$ and
$$
\U_{k_0}(\bar \mu, P,  \vert \bar Q\vert,r;\bar \nu ) \le 3^r \cdot (2k_0)^{4P+\vert \bar Q\vert }   \cdot 6^{3\mu_3}
\ \prod_{m=4}^{k_0} \, (2k_0)^{m\mu_m}  \ \prod_{k=k_0+1}^s (2k)^k.
\eqno (3.5)
$$

}

\vss
Lemma 3.1 is proved in Section 5 (see subsections 5.2 and 5.4).
Let us note that the condition (3.3) 
is the technical one related with the self-intersection degree of the root vertex (see subsection 5.1).

\section{Estimates of moments and proof of Theorem 2.2}

Let  $v^2= 1/4$ and $C_0$ be some constant specified later. 
Following the general  scheme of papers \cite{SS1,SS2,S} and \cite{R},
we represent the right-hand side of (2.10) as the sum of four  terms,
$$
 \E  \T (\hat A^{(n)})^{2s}  = \sum_{l=1}^4 Z_{2s}^{(l)},
$$
where
\begin{itemize}

\item $Z_{2s}^{(1)}$ is the sum over the subset  $\CI^{(1)}_{2s}$ of trajectories $ I_{2s}$ such that 
there is no multiple edge in the frame $\CF(w_{2s})$, $w_{2s} = w(I_{2s}), I_{2s} \in \CI_{2s}^{(1)}$ and  
$$
\vert \bar \nu (w_{2s}) \vert_1 = \sum_{k=2}^s (k-1)\nu_k \le C_0 s^2/n;
\eqno (4.1)
$$

\item $Z_{2s}^{(2)}$ is the sum over the subset  $\CI^{(2)}_{2s}$ of  $I_{2s}$ such that 
there exists at least one multiple edge in $\CF(w_{2s})$, $w_{2s} = w(I_{2s})$, condition (4.1) is verified, and 
 the maximal exit degree of $w_{2s}$  is bounded above by  $n^\delta$, \mbox{$D(w_{2s})\le n^{\delta}$},
 where  $\delta$ can be taken the same as in (2.5);

\item $Z_{2s}^{(3)}$ is the sum over the subset  $\CI^{(3)}_{2s}$ of  $I_{2s}$ such that 
there exists at least one multiple edge in $\CF(w_{2s})$, $w_{2s} = w(I_{2s})$, condition (4.1) is verified, and 
\mbox{$D(w_{2s})> n^{\delta}$;}

 \item $Z_{2s}^{(4)}$ is the sum over the subset  $\CI^{(4)}_{2s}$ of  $I_{2s}$ such that 
$\vert \bar \nu (w_{2s}) \vert_1 > C_0 s^2/n$, $w_{2s} = w(I_{2s})$.

\end{itemize}

Here and below we  denote $s= \lfloor s_n\rfloor$ and by even more abuse of denotations, 
we will write sometimes that $s=\chi n^{2/3}$. Also we will denote the limiting transition of (2.5)
as $(s,n)_\chi\to\infty$. 


\subsection{ Estimate of $Z^{(1)}_{2s}$}

As it is pointed out at the end of Section 3, one can restrict himself with the studies of the classes $\CC^{(\theta)}(\bar \nu, r)$ only.
Regarding a walk  $w_{2s}\in \CC^{(\th)}(\bar \nu,r)$, it is easy to see that 
the number of vertices in the graph $\Gamma'(w_{2s})$ is given by 
$$
\vert \CV(w_{2s})\vert = s+1 - \vert \bar \nu\vert_1.
\eqno (4.2)
$$
Then each class of trajectories  $\CI_{w_{2s}}$, $w_{2s}\in \CC^{(\theta)}(\bar \nu,r)$ is of the cardinality
$$
\vert \CI_{w_{2s}}\vert = n(n-1)\cdots (n-\vert \CV(w_{2s})\vert+1) = n(n-1)\cdots (n-s+ \vert \bar \nu\vert_1).
$$

We will use the following important estimate. 
\vs
\noindent {\bf Lemma 4.1. \cite{SS1} } {\it If $s<n$, then for arbitrary positive natural $\s<s$, the following inequality holds
$$
\prod_{k=1}^{s-\s} \left( 1 - {k\over n} \right) \le \exp\left\{ - {s^2\over 2n} \right\}\,  \exp\left\{ {s\s\over n}\right\}.
$$
}
\vs {\it Proof.}  One can use equality $1-k/n = \exp\{ \log (1-k/n)\}$ and then apply the Taylor expansion.

\vs

Regarding the sum over the classes of equivalence $\CC^{(\th)}(\bar \nu,r)$, using (3.2) together 
with the result  of Lemma 4.1 and
taking into account the obvious bound  for the weight $\hat \Pi(I_{2s})\le v^{2s}\le 4^{-s}$ (2.10), 
we write inequality
$$
Z^{(1)}_{2s} \le {1\over 4^s} 
\sum_{\th \in \Theta_{2s}} \sum_{\s=0}^{C_0s^2/n} 
\ \sum_{\bar \nu: \, \vert \bar \nu\vert_1 = \s}\ 
n{(n-1)\cdots (n-s+\s) \over n^{s-\s}} \cdot {1\over n^\s} \sum_{r=0}^{\nu_2}
\,  \vert \CC^{(\theta)}(\bar \nu,r)\vert  
$$
$$
\le 
n \exp\left\{ C_0{s^3\over n^2} - {s^2\over 2n}\right\} 
\ \sum_{\th\in \Theta_{2s}} \ 
\sum_{\s = 0}^{C_0s^2/n} \ 
\ \sum_{\bar \nu: \, \vert \bar \nu\vert_1 = \s}\ 
{1\over \nu_2!} \left( {s^2\over 2n} + {6sH_\th\over n}\right)^{\nu_2}
$$
$$
\times 
\prod_{k=3}^{s} {1\over \nu_k!} \left( { (2k)^k \, s^k \over k! \, n^{k-1} } \right)^{\nu_k}.
\eqno (4.3)
$$
Passing to the sum over all $\nu_i\ge 0, i\ge 2$ without any restriction, we obtain inequality
$$
Z_{2s}^{(1)} \le {n\over 4^s}  \, \exp\{ C_0 \chi^3\} \sum_{\th \in \Theta_{2s}} 
\exp\left\{ {6H_\th \chi^{3/2}\over \sqrt s} + 36\chi^3 
+ \sum_{k\ge 4} { (C_1 s)^k\over n^{k-1}}\right\},
$$
where $C_1 = \sup_{k\ge 2} {\displaystyle 2k \over \displaystyle (k!)^{1/k}}$. 
The Stirling formula implies that 
$$
n \vert \Theta_{2s}\vert = {n(2s)!\over s! \, (s+1)!} = {4^s\over \sqrt{\pi \chi^3}}(1+o(1)),
\quad (s,n)_\chi\to\infty.
$$
Then 
the following inequality holds (cf. (2.8))
$$
\limsup_{(s,n)_\chi \to\infty}  \ Z^{(1)}_{2s} \le\  {1\over \sqrt{\pi \chi^3}}\ B(6\chi^{3/2})\  e^{(C_0+36)\chi^3}.
\eqno (4.4)
$$
As a conclusion, it is easy to observe that the leading contribution to $Z_{2s}^{(1)}$ comes from the walks 
that have increasing number of 
simple self-intersections (open and closed ones) and triple self-intersections with no instants of open self-intersection at them.  
Also it is not hard to show that  the walks that have at least one loop edge provide vanishing contribution to 
$\limsup   \ Z^{(1)}_{2s}$.

The difference between the ensembles of real symmetric matrices and hermitian matrices is 
that the number of ways to pass the vertices open self-intersections by the non-marked exit steps is different for them
and in the case of hermitian matrices $B(6\chi^{3/2})$ can be replaced by $B(4\chi^{3/2})$. 


\subsection{Estimate of $Z^{(2)}_{2s}$}

Let  $k_0$ be an integer such that  
$$ 
k_0 \ge\,  {1\over  2\d}+1 > 4.
\eqno (4.5)
$$ 
Regarding a walk $w_{2s}$ that  belongs to a class 
$ \CC_{d,k_0}^{(\theta)}(\bar \mu, P,\bar Q,r ;\bar \nu^{(k_0)})$ 
determined in Section 3, we see that the condition to have at least one multiple  edge of $\Gamma'$ 
leads to inequality 
$$
\quad P + \vert \bar \nu^{(k_0)}\vert_1 \ge 1.
\eqno (4.6)
$$

It is easy to see that the weight of $I_{2s}$ such that 
$w(I_{2s})\in \CC_{d,k_0}^{(\theta)}(\bar \mu, P,\bar Q,r ;\bar \nu^{(k_0)})$ can be bounded as follows, 
$$
\hat \Pi(I_{2s}) \le\  {1\over 4^s}\cdot 
 (4V_4)^{P}\cdot   (2U_n)^{2\vert \bar Q\vert}\cdot  \prod_{k=k_0+1}^s \left(2U_n\right)^{2k\nu_k}\ .
\eqno (4.7)
$$
Finally, let us note that the number of vertices (4.2) is given now by relation
$$
\vert \CV(w_{2s})\vert  = s+1 - \sum_{m=1}^{k_0} (m-1) \mu_m - P - \vert \bar Q\vert - \sum_{k=k_0+1}^s (k-1)\nu_k.
$$

\vs

Multiplying the right-hand side of (3.4) by $n(n-1)\cdots (n-\vert \CV(w_{2s})\vert+1)/n^s$  
 and taking into account  (3.5) and  (4.7), 
we get the following upper bound for $Z^{(2)}_{2s}$,
$$
{Z^{(2)}_{2s} } \le {n\over 4^s}  
\ \sum_{\th\in \Theta_{2s}}\
\  \sum_{d=1}^{n^{\d}}
\ \  \sum_{\s = 0}^{C_0s^2/n} \ 
\  \sum_{
\stackrel{\bar \nu, \bar \mu,  P, \bar Q:  }
{  \vert \bar \mu\vert_1 +  P +\vert \bar Q\vert + \vert \bar \nu^{(k_0)}\vert_1 = \s} }^{(4.6)}\ 
\ \  e^{-s^2/n + C_0\chi^3}
$$
$$
 \times 
\
\ {1\over \mu_2!} 
\left( {s^2\over 2n} + {6sH_\th\over n} \right)^{\mu_2} 
\cdot {1\over \mu_3!} \left( {36s^3\over n^2}\right)^{\mu_3}\ 
\prod_{m=4}^{k_0} \, {1\over \mu_m!} \left( { (2k_0)^m s^m\over m! \, n^{m-1}} \right)^{\mu_m} 
$$
$$
\times \, {1\over P! } \left( { (4 k_0)^4 sdV_4\over n}\right)^{P} 
\cdot {1\over Q_1!} \left( { 16 k_0 P d U_n^2\over n}\right)^{Q_1}
$$
$$
\times\,  
 \prod_{j=2}^{k_0-2}  {1\over Q_j!} \left( {16k_0Q_{j-1} dU_n^2\over n} \right)^{Q_j}
 \ \   \prod_{k=k_0+1}^s\  
{1\over \nu_k!} \left( {(2C_1 s U_n^{2})^k\over  n^{k-1}}\right)^{\nu_k},
\eqno (4.8)
$$
where we have repeated some of the  computations  used in (4.3).

 \vs
 Let us point out that for all $n$ such that $n^{2\delta} > 16k_0 e \chi$,
 $$
 \sum_{Q_2, \dots , Q_{k_0-2} \ge 0}\    \prod_{j=2}^{k_0-2}  {1\over Q_j!} \left( {16k_0Q_{j-1} dU_n^2\over n} \right)^{Q_j}
\le e^{Q_1},
\eqno (4.9)
$$
thanks to (2.6) and the general bound $d \le s \le \chi n^{2/3}$.

\vs 
Remembering  the bound 
$
P\le  C_0\chi^2 n^{1/3}
$
we omit the condition  
$$
\vert \bar \mu\vert_1 +  P'  + \vert Q\vert + \vert \bar \nu^{(k_0)}\vert_1 = \s,
 \quad 0\le \s\le C_0s^2/n
$$
 and pass in the right-hand side of (4.8) to  the sum over all possible values of 
 $\bar \nu, \bar \mu,  P$ and $\bar Q$ such that (4.6) is verified. Taking into account (4.9)
 and using inequality $d\le n^\delta$, we get the bound 
 \begin{eqnarray*}
 Z^{(2)}_{2s} 
 &\le&{n^\d\over \sqrt{\pi \chi^3}}\, 
 B_s(6\chi^{3/2})  \, \exp\left\{ ( 36+C_0+ 36 )\chi^3 +\sum _{m=4}^{k_0} {(2k_0)^m \chi^m \over m!\,  n^{m/3-1}}
\right\} \\
 & &\times \left( \exp\left\{ {(4k_0)^4\chi e \over n^{1/3-\d}}V_4 
  +
  {1\over n^{2k_0\d - 1}}
 \sum_{k\ge k_0+1} 
 { (2C_1 \chi)^k\over n^{2(k-k_0)\d} }   \right\} - 1 \right),
 \end{eqnarray*}
 where the variable  $B_s$ is that of (2.9).
 Remembering definition   (4.5), we see that 
 $$
 Z^{(2)}_{2s} = o(1) \quad {\hbox {as}}\ \ (s,n)_\chi\to\infty.
 \eqno (4.10)
 $$


\subsection{Walks with high exit degree}

\subsubsection{Exit degrees and  proper and imported cells}

The sub-sum $Z^{(3)}_{2s}$ counts the contribution of the walks 
that have large maximal exit degree.
Let us take a particular walk $w_{2s}$ and denote by $\b_0$
the vertex whose exit degree is the  maximal one, i.e. such that $\vert \Delta(\b_0)\vert = D(w_{2s}) $.
Assume that $D(w_{2s}) = d$ and consider the edges of the exit cluster $\hat e_j \in \Delta(\b_0)$, $j=1,\dots, d$.

Regarding the Dyck structure $\theta_{2s} = \theta(w_{2s})$ and the  tree $T_s = T(\theta_{2s})$,
one can uniquely determine the set of edges $\Lambda = \{ \hat \vep_j, j=1,\dots, d\}$ of this tree $T_s$ that correspond  to the edges 
of the exit cluster $\Delta(\b_0)$. Each edge $\hat \vep_j\in \Lambda $ has a parent vertex $\hat \u_i$ of $T_s$ determined 
in obvious and natural way.  

\vsss
The edges of $\Lambda$ can have different parent vertices $\hat \u_i$ of $T_s$; we denote this number by $K$.  
For example, this can happen when  the self-intersection degree  $\kappa(\b_0)$ is greater than one,
or in other words, when there exists a number of marked arrivals at $\b_0$.
 In this case we say that  at $\b_0$ there exists a number of  {\it proper cells} (see subsection 5.6 for the rigorous definition). 
 
 Another possibility is given   by the case when $\b_0$ is visited
by a non-marked arrival instants of a special kind. 
To produce such non-marked arrivals at $\b_0$, the walk has to perform somewhere during its run 
an action 
that can be regarded as a break of the tree structure, or in other words, to have 
a non-zero number of broken tree structure instants of time, or briefly, the BTS-instants 
\cite{KV}.
In this case we say that $\b_0$ can have  a number of {\it imported cells}.

\vss
It is proved in \cite{KV} that the number $K$ of different parent vertices of the edges of $\Lambda$
is  bounded by the sum of $2\kappa (\b_0)$ and $L$, where $L$ is the total number of BTS-instants
performed by the walk during its run. We present the arguments of \cite{KV} in Section 5. 

\vss
Let us turn to the vertex of the maximal exit degree $\b_0$ and  assume  that 
the number of possible parent vertices of $\hat \vep_j$ is equal to  $K$.
It is easy to see that in this case the Dyck structure $\theta(w_{2s})$ is such that 
the corresponding tree $T_s = T(\theta_{2s})$ has at least one vertex
of the exit degree greater than $d/K$ \cite{SS2,S}. If the number $d/K$ is sufficiently large, then 
the set of  such Dyck paths is exponentially small with respect to the set of all possible Dyck paths $\Theta_{2s}$.
This makes possible to show that $Z^{(3)}$ vanishes in the limit $(s,n)_\chi\to\infty$.

\subsubsection{Placements of the BTS-instants}

The BTS-instants can be performed by the walk  at the instant of a self-intersection when the
walk arrives at the vertex that has at least two open edge attached to it, including the last arrival edge $e'$. 
If it is so, we will say that the BTS-instant is performed at the edge $e'$ (see Section 5 for more details).

Regarding the total number $L$ on BTS-instants performed by the walk,  
we can write that $L =  L_2 + L'$, 
where $L_2$ is the  number of BTS-instants performed by the edges at the vertices of self-intersection that have no
$p$-edges attached at them.
The remaining $L'$ BTS-instants can be produced at 
$p$-edges, 
$q$-edges, $\nu$-edges and 
by the  $\mu$-edges that are either edges with attributed $p$-edges or  the edges of non-simple $\mu$-self-intersections,.
Let us consider the cases of $\nu$-edges and $q$-edges first. 
\vs 
The number of BTS-instants  $L'_\nu$ that can be produced by the walk at the $\nu$-edges 
is such that $L'_\nu\le \sum_{k\ge k_0+1} (k-1)\nu_k = \vert \bar \nu^{(k_0)}\vert_1 $. 
Regarding the last two factors of (3.4) and (3.5) that  give the estimate of the number of walks 
with $\nu$-vertices, we can write that for any $h>1$
$$
\sum_{\bar \nu} \ {1\over h^{\vert \bar \nu^{(k_0)}\vert _1}} \   
\prod_{k=k_0+1}^s {1\over \nu_k!} \left( { (2C_1 s U_n^2)^k h^{k-1} \over n^{k-1}} \right)^{\nu_k} 
$$
$$
\le 
{1\over h^{L'_\nu}}\  \sum_{\bar \nu} \   
\prod_{k=k_0+1}^s {1\over \nu_k!} \left( { (2C_1 s U_n^2)^k h^{k-1}\over n^{k-1}} \right)^{\nu_k} . 
\eqno (4.11)
$$
\vs 
Regarding the BTS-instants $L'_q$ performed by the walk at the $q$-edges, we observe that
$
L'_q \le \vert \bar Q\vert \le (k_0-2) Q_1.
$ 
Therefore we can rewrite (4.9) in the following form with $h>1$,
$$
 \sum_{Q_1, \dots , Q_{k_0-2} \ge 0}\ {1\over h^{Q_1}}\cdot   {1\over Q_1!} \left( { d h P V_{6}\over n}\right)^{Q_1}
 \ \   \prod_{j=2}^{k_0-2}  {1\over Q_j!} \left( {16k_0Q_{j-1} dU_n^2\over n} \right)^{Q_j}
$$
$$
\le {1 \over h ^{L'_{q}/k_0}} \  e^{C_0 \chi^3 V_{6} h }.
\eqno (4.12)
$$
We  use  $V_{6}$ to simplify the arguments that will follow. In fact, 
the moments higher that $V_4$ can be avoided everywhere below but we do not take care about this. 

Having  (4.12),  we can also  include into consideration  the BTS-instants produced
by all  $\mu$-edges $\{\a_i,\b_i\}$ and attributed to them  $p$-edges that have $q$-edges over to  them. 
Indeed,
the number of such $\mu$- edges is bounded by $ Q_1$, the same concerns the corresponding $p$-edges, 
so the total number of edges of the form $\{\a_i,\b_I\}$ is bounded by $k_0 Q_1$.
Regarding the number $L'_{p,q}$ of BTS-instants produced by the edges of this form, we can write
inequality $L'_{p,q} \le k_0 Q_1$ and use (4.12) with $L'_q$ replaced by $L'_{p,q}$. That is why we have used
the factor $V_6$ in (4.12). 
From now on we can consider only the cases of $\mu$-edges that can have  $p$-edges 
attributed to them, but  have no $q$-edges over them.

\vs 
Let us consider the number of BTS-instants $L'_{p_1}$ that are produced by the \mbox{$p$-edges} attached to the $\mu$-edges
that arrive at the vertices of $\mu$-self-intersection degree 1. Then the corresponding self-intersections
are to be the open simple ones. Denoting the number of such edges by $P_1'$ and by $P_1''$ the number of 
other $p$-edges, we can write that 
$$
\sum_{P_1', P_1''}\, {1\over h^{P_1'}}  \cdot  {\mu_1!\over (\mu_1-P_1'-P_1'')! \, P_1'!\, P_1''!} 
\left( { 24h  H_\theta V_4\over n}\right)^{P_1'} \cdot \left( {4dV_4\over n}\right)^{P_1''} 
$$
$$
\le
{1\over h^{L'_{p_1} }} 
\exp\left\{ {4sdV_4\over n} + { 24hs  H_\theta V_4\over n}\right\}.  
\eqno (4.13)
$$
\vs 
Regarding the $\mu$-edges that participate at $\mu$-self-intersections of degree 2,  
we assume that there are $\mu_2'$ vertices whose $\mu$-edges 
have no $p$-edges attributed to them and there are $r_2$ vertices of open self-intersections there. 
The remaining $\mu_2''$ vertices of simple $\mu$-self-intersections are such that the corresponding
$\mu$-edges have $p$-edges attributed to them.

Then we can write that 
$$
\sum_{ r_2} 
{1\over h^{r_2}\, (\mu_2'-r_2)!} \left( {s^2\over 2n}\right)^{\mu_2' - r_2} \ 
{1\over r_2!} \left( { 6sh H_\theta\over n}\right)^{r_2} \le  
{1\over h^{L_2'}\, \mu_2'!} \, \left( {s^2\over 2n} + { 6sh H_\theta\over n}\right)^{\mu'_2},
\eqno (4.14)
$$
where $L_2'$ denotes the number of BTS-instants performed by the walk
at the instants of open simple self-intersections. 

Let us denote by  $L''_2= L'_{\mu''_2}$ the number of BTS-instants performed by the walk
at the edges of simple self-intersections and $p$-edges attached to them. 
It is clear that   $L''_2 \le 3\mu''_2$.  Then we can write the following inequality
$$
\sum_{\mu_2''} {1\over h^{\mu_2''}} \cdot {1\over \mu_2''!} \left( { 16h s^3 V_4\over n^2}
+{64 h s^4 V_4^2\over n^3}\right)^{\mu_2''} \le
{1\over h^{L''_{2}}} \ e^{ 16h \chi^3 V_4 + 64 h \chi^4 V_4^2 n^{-{1/3}}  }
\eqno (4.15)
$$
that is true for any  orientation of the $p$-edges.
\vss

Finally, for  $L_3'$ BTS-instants  performed by the walk at the $\mu$-edges attached to the vertices of $\kappa_\mu\ge 3$
and $P_3$ $p$-edges attributed to them, we can write the bound
$$
L'_3 \le P_3 + \sum_{m=3}^{k_0} (m-1)\mu_m.
$$
Then the sum over all such edges can be estimated as follows, 
$$
\sum_{P_3} \ \sum_{\mu_3, \dots , \mu_{k_0}} {3\mu_3 + \dots + k_0 \mu_{k_0} \choose P_3}\ 
\left( {4 dV_4\over n}\right)^{P_3}
\prod_{m=3} ^{k_0} \left( { (2k_0)^m s^m  \over m! \, n^{m-1}} \right)^{\mu_m} 
$$
$$
\le {1\over h^{L'_3} }\cdot  \exp
\left\{ \sum_{m=3}^{k_0} {4\chi V_4 m(2k_0)^m s^m h^{m-1}\over m! \, n^{m-2/3}} \right\}.
\eqno (4.16)
$$
Now we are ready to estimate the sub-sum $Z^{(3)}_{2s}$.

\subsubsection{Estimate of $Z^{(3)}_{2s}$}
Given integers $u$ and $d'$ from $ [1,\dots, s]$, we denote by $\Theta^{(u;d')}_{2s}$ the set of 
all Dyck paths $\theta$ such that $H_\theta= u$ and that $T(\theta)$ has at least one vertex 
of the exit degree not less than $d'$. 
Let us consider first the sub-sum $\tilde Z^{(3)}_{2s}$ over the subset of walks  such that  the vertex of the maximal exit degree
$\b_0$ is the $\mu$-vertex. Remembering that   $K\le 2k_0 +L$ and 
$L \le L_2' + L_2'' + L_3'+ L'_{p_1}   + L'_{q} + L'_{\nu}$ and  
taking into account relations (4.11)-(4.16),  we can write the following inequality
$$
\tilde Z^{(3)}_{2s} \le {n\over 4^s} \ e^{C_0\chi^3+h\chi^3 (16V_4  +  C_0 V_6) } \ \sum_{u=1}^s \ 
\ \sum_{d > n^\delta} \ \sum_{L\ge 0}\  \sum_{M=1}^{k_0}\ \  
\exp\left\{ {4sd V_4\over n}\right\}
$$
$$
\times {1\over h ^{2M+ L/k_0}}
\  \sum_{\theta \in \Theta^{(u;d/(2k_0+L))}_{2s}} \ \ \exp\left\{   { 6hs  H_\theta  ( 1 + 4 V_4)\over n} \right\}.
\eqno (4.17).
$$
It is proved in Section 5 that 
$$
\sum_{u=1}^s  \ \  \sum_{\theta \in \Theta^{(u;d')}_{2s} } \ \ \exp\left\{   { 6hs  H_\theta  ( 1 + 4 V_4)\over n} \right\}
\le 
6 s^2 \, e^{-\eta d'}\, t_s \, B_s( 6h \chi^{3/2} (1+ 4V_4)),
\eqno (4.18)
$$
where $\eta = \log (4/3)$ and $t_s =  \vert \Theta_{2s}\vert$.
Using (4.18), we get from (4.17) the upper bound
$$
\tilde Z^{(3)}_{2s} \le {n\over 4^s}\,  t_s \,  e^{C_0\chi^3+h\chi^3 (16 V_4  +  C_0 V_6) } \, B_s( 6h \chi^{3/2} (1+ 4V_4))
$$
$$
\times \, 6 s^3\,  k_0 \exp\left\{ {4sd V_4\over n}\right\}\  \sum_{L\ge 0} \exp\left\{ - { \eta d k_0\over 2+ L/k_0} - (2+L/k_0) \log h\right\}.
\eqno (4.19)
$$
Elementary analysis of the function
$$
f(y) =  { \eta d k_0\over 2+ y/k_0} +  (2+y/k_0) \log h, \quad y>0
\eqno (4.20)
$$
shows that its minimum value is given by $2\sqrt{ \eta d k_0\, \log h}$. Therefore
$$
\tilde Z^{(3)}_{2s} \le {6e^{C_2}k_0\over \sqrt{\pi \chi^3}} \,  B_s( 6h \chi^{3/2} (1+ 4V_4)) \, s^4 
e^{ - 2n^{\delta/2} \left( \sqrt {\eta k_0 \, \log h}  - {2\chi^{3/2}  V_4 } \right) },
\eqno (4.21)
$$
where $C_2 = C_0\chi^3+h\chi^3 (16 V_4  +  C_0 V_6)$. 

Taking $h = h_0+1$ such that $\log h_0 =  4\chi^3 V_4^2/(\eta k_0)$, we conclude that 
$$
\tilde Z^{(3)}_{2s} = o(1), \quad (s,n)_\chi\to\infty.
\eqno (4.22)
 $$
The first part of the estimate of $Z^{(3)}_{2s}$ is completed.

The second part of the reasoning concerns the sub-sum $\hat Z^{(3)}_{2s}$ over the walks such that 
the vertex of maximal $\b_0$ exit degree
is the $\nu$-vertex of the self-intersection degree $N\ge k_0+1$. In this case the proof of the relation 
$$
\hat Z^{(3)}_{2s} = o(1), \quad (s,n)_\chi\to\infty
\eqno (4.23)
$$
can be obtained by almost the same arguments as those used to estimate $\tilde Z^{(3)}_{2s}$. 
Indeed, the upper bound for $\hat Z_{2s}^{(3)}$ is similar to (4.17) and is of the form
$$
\hat  Z^{(3)}_{2s} \le {n\over 4^s} \ e^{C_0\chi^3+h\chi^3 (16V_4  +  C_0 V_6) } \ \sum_{u=1}^s \ 
\ \sum_{d > n^\delta} \ \sum_{L\ge 0}\  \sum_{N=k_0+1}^{s}\ \  
\exp\left\{ {4sd V_4\over n}\right\}
$$
$$
\times {1\over h ^{2N+ L/k_0}} \cdot {(2C_1 s h^2 U_n^2)^N \over n^{N-1} } \  
\  \sum_{\theta \in \Theta^{(u;d/(2N+L))}_{2s}} \ \ \exp\left\{   { 6hs  H_\theta  ( 1 + 4 V_4)\over n} \right\}.
\eqno (4.24)
$$
It is easy to see that after the use of (4.18), we arrive at the analysis of the function of the form (4.19)
where $2+ y/k_0$ is replaced by $z = 2N + L/k_0$. Then we obtain from (4.24) the estimate of the form (4.20),
 where $k_0 s^4$ is replaced by $s^5$ and (4.23) follows. 
We do not present 
the details of computations.


\subsection{ Estimate of $Z^{(4)}_{2s}$}

In this subsection, we mainly  follow the lines of the proof indicated in \cite{R} and 
completed in  \cite{KV}. We do not need here to determine the  $\mu$-structure of  
$w_{2s}$. 

\vs
\noindent {\bf Lemma 4.1.} 
{\it Given any walk $w_{2s}$ of the class 
$\CC(\bar \nu)$, the weight $\hat  \Pi(w_{2s}) $  (2.10) is bounded as follows;
 $$
\hat  \Pi(w_{2s}) \le  {1\over 4^s} \cdot \prod_{k=2}^s \left(16V_{12}\cdot   (2U_n)^{2(k-2)}
 \right)^{\nu_k}.
 \eqno (4.25)
  $$}
  
 {\it Proof.} Regarding  a vertex $\g$ with
 $\kappa(\g)\ge 2 $ of 
  $\Gamma(w_{2s}) = (\CV,\CE)$,
 we color in red the first two marked
arrival edges at $\g$ and their non-marked closures. Passing to
another vertex with $\kappa\ge 2$, we repeat the same procedure
and finally get $4\sum_{k=2}^s\nu_k$ red edges. 
We also color in red the first marked arrival edge at   the root vertex $\rho$,
as well as its non-marked closure. 
Regarding the vertices $\b$ with $\kappa(\b)=1$, we color in yellow the marked arrival edges and
their non-marked counterparts.

In $\Gamma(w_{2s})$, it
remains $2\sum_{k=2}^s (k-2)\nu_k$ non-colored (grey) edges. 
Regarding the weight $\hat \Pi(w_{2s})$, let us replace the random variables
that correspond to the grey edges by their upper bounds $U_n$ (2.6).
Then we get inequality
$$
\hat \Pi(w_{2s}) \le {1\over 4^s} \ (2U_n)^{2(k-2)\nu_k}\cdot \hat  \Pi^* (w_{2s}),
$$
where $\hat \Pi^*(w_{2s})$ represents the product of the
mathematical expectations of the random variables associated with
the colored edges of $\CE(w_{2s})$ normalized by a suitable number of factors $V_2^{-1}$.
Let us denote by $\Gamma^* = (\CV,\CE^*)$ the graph with these oriented marked colored edges. 

\vsss
Let $\tilde \b$ be a vertex of $\Gamma^*$ that serves as the end of two red edges we denote by  
$(\a',\tilde \b)$ and $(\a'',\tilde  \b)$. Let us first assume that $\a'\neq \a''$. Then it is not hard to see
that there can be only one red or blue edge $(\tilde  \b,\a')\in \CE^*$. If it would be two red edges 
$(\tilde \b,\a')$,
then $\a'$ should be the root vertex $\rho$, but we assumed from the very beginning
 that $\rho$ always has already one invisible
arrival edge;  therefore  one of these two red edges is  actually the grey one absent in $\CE^*$. 
The same concerns the vertex $\a''$. Thus in this case the weight contribution of 
the edges that end at $\tilde  \b$  to $\hat \Pi^*(w_{2s})$  is bounded by $16V_4^2\le V_6^2$.  We assume
without loss of generality  that 
$1/4=V_2<1\le V_4\le V_6\le V_{12}$.

The same reasoning shows that if $\a'=\a''=\a$, then  it  exists  only one red edge  of the form $(\tilde \b,\a)$.
Then the edges that end at $\tilde \b$ provide the contribution not greater than $V_6$. 
The total number of such vertices $\tilde \b$ is bounded by $2\sum_{k\ge 2} \nu_k$. Then clearly 
$\hat \Pi^* (w_{2s}) \le \prod_{k\ge 2} (16V_{12})^{\nu_k}$ and (4.25) follows. Lemma 4.1 is proved. 

\vss

Let us show that
$$
Z^{(4)}_{2s} = o(1), \quad (s,n)_\chi\to\infty.
\eqno (4.26)
$$ 
Using  inequality (3.1) and taking into account   (4.25), we can write  that
$$
Z^{(4)}_{2s} \le {(2s)!\over 4^s\, s!\, (s+1)!} 
\ \sum_{\s\ge C_0s^2/n}
\ \  \sum_{\bar \nu: \, \vert \bar \nu\vert_1 = \s}\ 
{ n(n-1)\cdots (n-s+\s)\over n^s}
$$
$$
\times\ 
\prod_{k=2}^s {s^{k\nu_k} \over \nu_k!} \left(4C_1U_n^2\right)^{(k-2)\nu_k} \left( 16C_1^2V_{12}\right)^{\nu_k}.
\eqno (4.27)
$$
Denoting $\vert \bar \nu\vert_2 = \sum_{k=2}^s (k-2)\nu_k$, using identity 
$n^s = n^{s-\s} n^\s$  and taking into account that 
$$
\prod _{k=2}^{s} \, s^{k\nu_k} = s^{2 \vert \bar \nu\vert_1 - \vert \bar \nu\vert_2}\ ,
$$
we obtain from (4.27) that
$$
Z^{(4)}_{2s} \le  {n(2s)!\over 4^s\, s!\, (s+1)!} \ \sum_{\s\ge C_0s^2/n}
\ \sum_{\z =0}^\s \ {1\over n^\s} \cdot {s^{2\s}\over s^\z}  \left( 16C_1^2V_{12}\right)^{\s-\z}
$$
$$
\times \ \  \sum_{{\bar \nu: }{ \vert \bar \nu\vert_1 = \s, \vert \bar \nu\vert_2 = \z}}\ \ \  
\prod_{k=2}^s {1 \over \nu_k!} \left(4C_1U_n^2\right)^{(k-2)\nu_k}.
\eqno (4.28)
$$ 
Multiplying and dividing the right-hand side of (4.28) by $\s! /(\s-\z)!$, we obtain that 
$$
Z^{(4)}_{2s}\le  {n(2s)!\over 4^s\, s!\, (s+1)!} \
\sum_{\s\ge C_0s^2/n} \ {1\over \s!} \left( { s^2\over n}\right)^{\s} 
\sum_{\zeta=0}^\s 
 {\s!\over (\s-\z)!\, \s ^\z} \cdot \left( 16C_1^2 V_{12}\right)^{\s-\z}
$$
$$
\times \ 
\  \sum_{\bar \nu: \ \vert \bar \nu\vert_1 = \s, \ \vert \bar \nu\vert_2 = \zeta} \ \ 
{ (\s-\zeta)! \over \nu_2!\cdots \nu_s!}  
\ \  \prod_{k\ge 2} \left({16C_1 \s U_n^2\over s}\right)^{(k-2)\nu_k}. 
\eqno (4.29)
$$

\vs
Now we separate the sum over $\s$ into two parts:
$$
\CS_1 = \left\{ \s: C_0\chi^2 n^{1/3}\le \s\le {\chi n^{1/3+2\d}\over 32C_1}
\right\},
$$
$$
\CS_2 = \left\{ \s:  {\chi n^{1/3+2\d}\over 32C_1} < \s \le s\right\}
$$
and denote by $Z^{(4;i)}_{2s}$, $i=1,2$ the corresponding sub-sums.

\vskip 0.4cm

If $\s\in \CS_1$, then 
$$
  \sum_{\ \vert \bar \nu\vert_1 = \s, \ \vert \bar \nu\vert_2 = \zeta} \ \ 
{ (\s-\zeta)! \over \nu_2!\cdots \nu_s!}  
\ \  \prod_{k\ge 2} \left({16C_1 \s \over \chi n^{1/3+2\d} }\right)^{(k-2)\nu_k}
$$
$$
\le \sum_{\nu_2+\dots + \nu_s = \s-\zeta } \ 
{ (\s-\zeta)! \over \nu_2!\cdots \nu_s!}\   \prod_{k=2}^s {1\over 2^{(k-2)\nu_k}} \le 2^{\s-\zeta},
$$
where we have used the multinomial theorem.  
Then we can  derive from (4.29) the bound
$$
Z^{(4;1)}_{2s}\le
 {n(2s)!\over 4^s\,  s!\, (s+1)!}\ \ 
\sum_{\s\in \CS_1} \ {1\over \s!} \left( { s^2\over n}\right)^{\s}
\ \ \sum_{\zeta=0}^\s  
{\s!\over (\s-\z)!\,  \z!}\,  (32C_1^2 V_{12})^{\s-\z}\, .
$$
\vs
\noindent Using inequality  $\s! \ge (\s/e)^\s \sqrt{2\pi \s}$, we conclude  that 
$$
Z^{(4;1)}_{2s}\le
{n(2s)!\over 4^s \,  s!\, (s+1)!} \cdot {1\over \sqrt {2\pi C_0 \chi^2 n^{1/3}}}\ 
 \sum_{\s\in \CS_1} 
 \left( {s^2 e\over n \s} \left( 1+32C_1^2 V_{12}\right) \right)^\s.
 $$
Remembering that  $\s\ge C_0 s^2/n$, we see that the last series converges provided
$$
C_0> e(1+ 32C_1^2\,  V_{12}) 
\eqno (4.30)
$$
and then $Z^{(4;1)}_{2s} = o(1)$ in the limit $(s,n)_\chi\to \infty$. 

\vskip 0.4cm

Let us consider $Z^{(4;2)}_{2s}$. It follows from (4.29) that 
$$
Z^{(4;2)}_{2s} \le { n(1+o(1)) \over \sqrt{\pi s^3}} \sum_{\s\in \CS_2} {1\over \s!}
\left( { 16s^2 C_1 V_{12} \over n} \right)^{\s} \ \sum_{\z = 0}^\s \left( { U_n^2\over C_1 V_{12} s} \right)^\z 
$$
$$
\times {\s!\over (\s-\z)!} \  \  \sum_{\bar \nu: \ \vert \bar \nu\vert_1 = \s, \ \vert \bar \nu\vert_2 = \zeta} \ \ 
{ (\s-\zeta)! \over \nu_2!\cdots \nu_s!}  .
\eqno (4.31)
$$
It is not hard to see that for any real $b_i$ the following holds,
$$
\  \sum_{\bar \nu: \ \vert \bar \nu\vert_1 = \s, \ \vert \bar \nu\vert_2 = \zeta} \ \ 
{ (\s-\zeta)! \over \nu_2!\cdots \nu_s!}\  b_2^{\nu_2}\cdots b_s^{\nu_s} 
$$
$$
\le 
\ 2^\s \,  \sum_{\nu_2+\dots +\nu_s = \s-\z} \ \ 
{ (\s-\zeta)! \over \nu_2!\cdots \nu_s!}\  \left({b_2\over 2}\right)^{\nu_2}\cdots  \left({b_s\over 2^{s-1}}\right)^{\nu_s}
$$
$$
= 
2^\s \left( {b_2\over 2} + \dots + { b_s\over 2^{s-1}}\right)^{\s-\z}.
\eqno (4.32)
$$

Using (4.32) with $b_i=1$, we get from (4.31) the bound
$$
Z^{(4;1)}_{2s} \le  { n(1+o(1)) \over \sqrt{\pi s^3}} \sum_{\s\in \CS_2} \, 
\left( { 32C_1V_{12} e s^2\over n \s} \right)^\s \,
$$
$$
\times  {2\over \sqrt{ 2\pi \s}} \ \ \sum_{\z =0}^\s 
\left( { \s\over C_1 e V_{12} \chi n^{1/3+2\d}}\right)^\z { \s! \sqrt{2\pi \z} \over (\s-\z)! \, \z!},
$$
where we have used inequality $\z! \le \sqrt{8\pi\z} (\z/e)^\z$. 
Then we can write that
$$
Z^{(4;2)}_{2s} \le {1+o(1)\over \sqrt{\pi \chi^3}} \ \sum_{\s \in \CS_2} 
\left( { 64 \chi\over n^{\d}} \right)^\s \sum_{\z'=0}^\s 
\left( { C_1 e V_{12} \chi n^{1/3+\d}\over \s} \right)^{\z'}.
$$
The last series converges and then  $Z^{(4;2)}_{2s} = o(1)$ in the limit $(s,n)_\chi\to \infty$. 
This fact completes the proof of  (4.26).

\subsection{Proof of Theorem 2.2}

Regarding  (4.4) together with relations (4.10), (4.22), (4.23)  and (4.26), we conclude that  inequality (2.7)
 is true.  
Returning to the moment 
$\hat M_{2s}^{(n)}=  \E \T ( \hat A^{(n)})^{2s}$, we see that the leading contribution to 
$\hat M_{2s}^{(n)}$ 
is provided by the sub-sum
$Z^{(1)}_{2s}$ over the walks that have no multiple edges; 
the sum over subset of walks with multiple edges  given by $Z^{(2)}_{2s}+ Z^{(3)}_{2s}+ Z^{(4)}_{2s}$
 vanishes in the limit $s_n = \chi n^{2/3}, n\to\infty$. 
Also the walks whose graphs have at least one loop do not contribute
to the limiting expression of
 $\hat M_{2s}^{(n)}$.
Therefore one can
conclude that although the upper bound (2.7) depends on the technical constant 
$C_0$ (4.30),
the leading contribution to 
$
\hat M_{2\lfloor s_n\rfloor }^{(n)} $ does not depend
on particular values of the moments $V_{2k} = \E a_{ij}^{2k}$. 
This implies the universality of the upper bound 
$
  \limsup_{s_n= \chi n^{2/3},\,  n\to\infty}\E \T (  \hat A^{(n)})^{2\lfloor s_n\rfloor }.
  $

Moreover, the leading contribution to the moments of Wigner random matrices 
$\hat M_{2s}^{(n)}$ coincides with that obtained for the moments of 
the matrices of Gaussian Orthogonal Ensemble. In the latter case, the limit 
of $M_{2s_n}^{(n)}$ (2.4) is known to exist \cite{S}. Then the limit (2.6) also exists.
Theorem 2.2 is proved.

\vss 

{\it Remark.} Here we are in the situation similar to that  encountered in paper \cite{S},
where the upper bound of the form of (2.8) is obtained under condition
(1.4).
 While  it is not clearly indicated in \cite{S},
the upper bounds obtained there  depend on the value of the constant $C_a$.
However, the leading contribution to the corresponding moments  $ M_{2s}^{(n)}$ 
comes from the walks with no multiple edges. 
This  implies  universality of $\limsup_{(s,n)\chi\to\infty}  M_{2s}^{(n)}$ with respect 
to the probability distribution of random variables $a_{ij}$.
The same concern the estimates given in \cite{R} that should involve the constant $C_b$ of (1.5).


\subsection{Proof of Theorem 2.1}

Using the standard arguments of the probability theory,
we can write that 
$$
 \P \left\{ \hat \l_{\max}^{(n)} > 2v \left( 1+ {x\over n^{2/3}}\right) \right\} \le
  { \E \T ( \hat A^{(n)})^{2s_n} \over 
 \left( 2v(1+x n^{-2/3})\right)^{2s_n}},
 \eqno (4.33)
 $$
 where $\hat \l_{\max}^{(n)} = \l_{\max}(\hat A^{(n)})$.
 Then, regarding the limit $s_n,n\to\infty$, $ s_n=\chi n^{2/3}$, we deduce with the help of (2.7)
 the following bound 
 $$
 \limsup_{n\to\infty}  \P \left\{ \hat \l_{\max}^{(n)} > 2v \left( 1+ {x\over n^{2/3}}\right) \right\} 
 \le \inf_{\chi>0} {\cal G}(\chi ) e^{-x\chi},
 \eqno (4.34)
 $$
 where 
 ${\cal G} (\chi ) = B(6\chi^{2/3} ) e^{C\chi^3} ( \pi \chi^3)^{-1/2}$. Due to  Theorem 2.2,
 the upper bound (4.34) is also true with ${\cal G}(\chi) $ replaced by $\CL(\chi)$. 
 
\vss
Let us consider the subset 
${\cal O} _n
= \cap_{1\le i\le j\le n} \left\{ \omega:  \vert a_{ij}\vert \le n^{1/6 - \delta}\right\} \subseteq \Omega$.
Then  
$$
\P \left\{ \l_{\max}(A^{(n)}) > y\right\} = 
\P \left\{ \hat \l_{\max}^{(n)} > y\right\}  + \P \left\{ \left(\l_{\max}(A^{(n)}) > y\right) \cap \overline {\cal{ O}}_n\right\} .
\eqno (4.35)
$$
It is easy to prove that 
$$
\P \left( \overline  {\cal O}_n \right) = o(1), \quad n\to\infty
\eqno (4.36)
$$
under condition (2.3)  of Theorem 2.1 (see  subsection 5.7). Then  (4.35)
combined with inequality (4.33) and convergence (2.7) implies the bound (2.4). Theorem 2.1 is proved.

\section{Background statements and auxiliary estimates}

\subsection {Formula (3.1) for the root vertex }

Let us start with a version of the formula (3.1) that takes into account the $N$ marked arrival instants at the root vertex $\rho$.
In this case the number of possibilities to produce a partition of the set $\Xi_s$ is given by
$$
\Psi_s(\bar \nu, N) = {s!\over (s-N-\sum_{k\ge 2} k\nu_k)! \cdot N! \cdot \prod_{k\ge 2} (k!)^{\nu_k}\, \nu_k!} \ .
$$
Regarding the $\nu_{N+1}$ vertices of self-intersection of the order $k'=N+1$, we observe that 
$$
{s'!\over (s'-N-k'\nu_{k'})! \cdot N! \cdot (k'!)^{\nu_{k'}} \, \nu_{k'}!  }\le {s'!\over (k'!)^{\nu_{k'}+1} \, (\nu_{k'}+1)!}
$$
provided $2(k'+1)\nu_{k'} \le s' -1= s-\sum_{k\ge 2, k\neq k'}k\nu_k$ (cf. (3.3)). This condition is fulfilled in the case (4.2) we consider.
Therefore, the left-hand side of  formula (3.1) is in agreement with the definition of the self-intersection degree 
of the root vertex and can be used without restrictions.

\subsection{Proof of  (3.5) }

\n {\bf Lemma 5.1.}  {\it Consider a walk $w_{2s}$ whose  graph $\Gamma'(w_{2s})$ contains
a vertex $\b$ that has $m$,
 $p_\b$ and $q_\b$ 
arrival $\mu$-edges,   $p$-edges and  $q$-edges, respectively. 
 The following propositions are true:
\vss 
\noindent a) the total number of non-marked edges of the form $(\b,\g_i)$ is equal to $m+p_\b+q_\b$;
\vss
\noindent b) at any instant of time $t$ 
the  number of marked edges attached to  $\b$ is bounded by $2\kappa(\b)$.
}
\vs
{\it Proof.} We prove Lemma 5.1 by using the reasoning of \cite{KV}. 
Let us  denote by $I$ and $\bar I$ the numbers of marked and non-marked edges that arrive at $\b$ and by 
$O$ and $\bar O$ the numbers of marked and non-marked edges that leave $\b$. The walk $w_{2s}$ is closed and even,
and this implies relations
$$
I+ \bar I = O + \bar O, \quad {\hbox{and }} \quad I+O = \bar I + \bar O
$$
that lead to equality 
$\bar I - O = O - \bar I$. Then $ \bar I = O$ and   
$$
\bar O = I = m+p'_\b+q_\b.
\eqno (5.1)
$$ 
This proves the part (a) of Lemma 5.1. 

We  prove the second statement of Lemma 5.2 for the case when $\Gamma(w_{2s})$ has no loops, i.e. 
the edges of the form $(\b,\b)$. 
Let us introduce the variables
$I(t)$, $\bar I(t)$, $O(t)$, and $\bar O(t)$ that count the corresponding number of marked and non-marked arrivals
to $\b$ 
and marked and non-marked departures  from $\b$ during the time interval $[0,t-1]$.
Regarding the number of the edges $A_\b(t)$ attached to $\b$ that are $t$-open in the sense of Definition 2.1, one 
gets equality
$$
A_\b(t) = I(t)+ O(t)  - \bar I(t) - \bar O(t) .
\eqno (5.2)
$$
If the walk has left $\b$ at the instant $t$, then
$
I(t) + \bar I(t) = O(t) + \bar O(t) - \phi,
$
where $\phi$ equals to $0$ or $1$ in dependence whether $\b$ is the root vertex of $\Gamma(w_{2s})$ or not. 
Then we get equality $O(t) = I(t) + \bar I(t) - \bar O(t) + \phi$ that together with (5.2) implies  relations 
$$
A_\b(t)   = 2 I(t)  - 2\bar O(t) + \phi \le 2I(t)+\phi.
\eqno (5.3)
$$

If $\b$ is not the root vertex, then $\phi =0$ and $I(t)\le I(2s)= 2\kappa(\b)$. 
If $\b=\rho$ is the root vertex, then $\phi =1$. However, since we assume the presence of one invisible
arrival edge at the root vertex $\rho$, then $\kappa(\rho) = I(2s)+1$ and (5.3) again implies inequality 
 $A_\b(t) \le 2\kappa(\b)$. 

\vs Using the same argument, it is easy to show that 
the statement (b) remains valid also in the case when $\Gamma(w_{2s})$ has loops. Lemma 5.1 is proved.

{\it Remark.} It is not hard to show that in the case of simple self-intersection the number of open edges attached
to the vertex $\b$ at any arrival instant $t$ is equal to 1 in the case when the self-intersection is not open.
If the self-intersection is the open one, than the number of the open edges 
attached at $\b$ at the time of the second arrival by a marked edge is equal to 3. Then 
a number of repetitions of pairs of steps $(O,\bar I)$ can happen and at the instant of the step
$\bar O$ the walk has 3 possible ways to go. Then again a number of couples $(O,\bar I)$ can take place
and at the instant of the second non-marked departure the walk has only one possibility to close the remaining open 
edge attached at this vertex. 

\vs 
{\it Proof of  inequality (3.5).} 
It follows from Lemma 5.1 that for any vertex $\b$ of the graph of the walk $w_{2s}\in \CC_{d,k_0}(\bar \mu, P, Q,r)$,
such that $\kappa_{\mu}(\b)\ge 4$, we have
$$
\Upsilon(\beta) \le (2k_0)^{m+p_\b+q_\b}.
\eqno (5.4)
$$
It is not hard to see that if $\kappa_\mu(\b)=1$, then the non-trivial case is given by $p_\b=1$, 
and we can write that
$
\Upsilon(\b)\le  (2k_0)^{2+ q_\b}.
$

The remaining  two cases with  $p_\b\ge 1$ are described by  the relations
$$
\Upsilon(\b)\le  (2k_0)^{3+ q_\b}, \quad \hbox{if \ } \kappa_\mu(\b)=2 
\eqno (5.5)
$$
and 
$$
\Upsilon(\b)\le  (2k_0)^{4+ q_\b}, \quad \hbox{if \ } \kappa_\mu(\b)=3. 
\eqno (5.6)
$$
If $\kappa_\mu(\b)=3$ and $p_\b'=0$, then $\Upsilon(\b)\le 6$.

Regarding the product over  all vertices $\b\in {\cal E}_g$ 
that takes into account all $p$- and $q$-edges together with their orientations,
and taking into account the fact that 
simple open $\nu$-intersections give factors $3$, we get formula (3.5).

\subsection{Multiple edges and the last passage principle }

Let us  repeat important  elements of reasoning used in \cite{S} to  prove the upper bound  (3.2).
We assume   for simplicity  that  $\nu_3 = \dots = \nu_s =0$ and consider first the simplest case
of $r=p=1$. 
Given $\theta\in \Theta_{2s} $, let us consider the family of walks
such that  the numbers $0<x_1< \dots < x_{\nu_2}<s$ determine the instants 
$1\le \xi_{x_j}\le 2s-1$ of simple self-intersections. 

Among these $\mu_2$ values, we choose 
two values, $y$ and $z$ that will be the instants of the open simple self-intersection and
of the self-intersection that produces the multiple edge, respectively and assume that $y<z$. 
For each of the  remaining $\nu_2 - 2$ values $x_j$, we choose the value $v_j<x_j$
that determines the vertex of the simple self-intersection. 

Having determined these parameters, we start the run of the walk according to $\theta$
and the values prescribed until the instant of time
$\xi_y-1$. At this instant of time the walk, if it exists, staying at the vertex $\a= w(\xi_y-1)$ has to join
a vertex $\b$ such that there exists at least one open edge attached to $\b$. 
The principal observation here is that there is not more than $2 \theta(\xi_y-1)$ vertices of this kind \cite{SS2,S}.
Then the number of choices of the vertex $v_y$ of the open self-intersection
is bounded by $ 2H_\theta$, where $H_\theta$ is the maximal level reached by the Dyck path $\theta$. 
 
 When the vertex $v_y$ is chosen, the walk draws the edge $(\a,v_y)\in \Gamma$.
 If the next step $\xi_y+1$ is non-marked, the walk has to decide where to go  according to the rule 
 prescribed by $\Upsilon$.  If the step $\xi_y+1$ is marked, the walk continues its run until it returns to the 
 vertex $v_y$ with the non-marked step out of it. 

\vs 
The walk continues it run till the instant of time $\xi_z-1$ when it arrives at the vertex $\g= w(\xi_z-1)$.
At the marked instant of time $\xi_z$, it looks to join a vertex $\b'$ such that the closed edge $\{\g,\b'\}$
exists. This can be either a vertex of the exit cluster $\Delta_\g$ of the vertex of the marked arrival
$(\b'',\g)$. In any case, the part $[0,\xi_z-1]$ of the walk being completely determined by the choice
of the vertex $v_y$ that by the rule $\Upsilon$, the set $\Delta _\g$ is determined and has the cardinality
not greater than $d$. Also the marked arrival edges of the form $(\b'',\g)$ are completely determined; 
their number is bounded by $2$ in our case. If the self-intersection at the instant $z$ is prescribed
to be open, then the number of choices of $v_z$ is bounded by $\tilde d= \max\{ 2H_\theta, d+2\}$. 
Then the walk continues its run (again, if it is possible) till the final instant of time $2s$. We see that given the values
of the $\nu_2-2$ self-intersections, there is not more than $2H_\theta \times \tilde d$  possible walks. 
\vs 
Using this base observation of \cite{SS2,S}, it is not hard to prove (3.2) by recurrence with respect to  $r,p\ge 1$ and to consider
the general case of $\bar \nu$. 
Important thing is that in this reasoning, the last arrivals at the vertices of open simple self-intersections
and simple self-intersections producing  double edges are to be fixed. 

\vs The natural generalization of this rule of the last arrival to the case of multiple edges is given by the 
principle of the last passage. We will need it in the proof of estimate (3.4). 
Let us consider a family of walks  $w_{2s} \in \CC^{(\theta)}_d(\bar \nu)$  such that the graph $\Gamma'(w_{2s})$ with marked edges  contains a couple 
$\{\a,\b\}$ such that there are $g>1$ marked edges joining them. 
We consider the walks whose vertices $\g$ are such that $\kappa(\g)\le k_0$.

We choose the value $y$, $1<y\le s$ 
of the last passage of this couple by a marked step $1<\xi_y\le 2s-1$  and start the run of the walk
on the time  interval $[0,\xi_y]$. We do this 
according to the part of self-intersections that is prescribed to be performed during $[0,\xi_y]$
and the rules $\Upsilon$ of choices at the vertices of the open self-intersections.  

Assuming that $w(\xi_y-1) = \a$ and $w(\xi_y)=\b$, we ask about the number of values that can be seen
at remaining $l-1$ marked edges $\{\a,\b\}$ during the run of the walk. Since the part 
$[0, \xi_y]$ of $w$ is determined, the exit clusters $\Delta_\a(\xi_y)$ and $\Delta_\b(\xi_y)$
are determined, as well as the marked closed edges of the form
$(\g',\a)$ and  $(\g'',\b)$. Then there are not more than $d$ or $k_0$ values to choose in dependence of the orientation
of one or another of $l-1$ edges. 

One can repeat  this reasoning in  the case of the walks with  a number of different multiple edges.
Assuming that the instants of the last passages $\xi_{y_1}<\dots < \xi_{y_k}$ as well as the instants of  all other self-intersections
are communicated to us, 
someone produces a realization of a walk 
from the given class by keeping  shadowed from us the values of  the remaining $l_1-1, \dots , l_k-1$ passages.
Then, starting from the last multiple edge $e(\xi_{y_k})$, we estimate step by step the number of possible values.

\subsection{Proof of Lemma 3.1}

Let us consider the class of walks $\CC^{(\theta)}_{d,k_0}(\bar \mu, P, \bar Q,r;\bar \nu)$. 
We can order the $\mu$-edges $(\a,\b)$ of $\Gamma'(w_{2s})$ according to the order of the
last passages of the corresponding pairs $\{\a,\b\}$ that will be created after that 
$p$-edges and $q$-edges will be distributed over the $\mu$-edges.

We consider 
 $\vert \bar \mu\vert' = \mu_1+ (\mu_2-r) + \sum_{m=3}^{k_0} m\mu_m$
ordered $\mu$-edges and choose $P$ to put $p$-edges over them. 
The number of choices is estimated by the following inequalities, 
$$
2^P { \vert \bar \mu \vert' \choose P} \le { (2\vert \bar \mu'\vert) ^P\over P!} \le {(2 s)^P\over P!},
\eqno (5.7)
$$
where $2^P$ takes into account the orientations of $p$-edges. 
The next step is to put  $Q_1$ edges  over the $P$ green edges. This can be done
in not more than 
$$
2^{Q_1} { P \choose Q_1}  \le {(2 P)^{Q_1}\over Q_1!}
\eqno (5.8)
$$
ways, where $2^{Q_1}$ takes into account the orientations. 
We continue these steps till the last $Q_{k_0-2}$ edges
are placed. 

Having all $p$- and $q$-edges distributed, we choose the instants of the last passages of the edges attributed to the $\mu$-edges.
According to the combinatorial reasoning of (3.1), this can be done by not more than
$$
{1\over (\mu_2-r)! } \left( {s^2\over 2} \right)^{\mu_2-r} \cdot {s^r\over r!} \ \prod_{m=3}^{k_0} 
{1\over  \mu_m!} \left( { s^m\over m!} \right)^{\mu_m}
\eqno (5.9)
$$
ways, where the numbers of values of the last arrivals at the vertices of open simple self-intersections are also estimated.
Similar expression estimates the number of possible values at the vertices of $k$-fold self-intersections with $k\ge k_0+1$.

The principle of last passage instants 
says that the number of possible values at the edges lying over the multiple edges is bounded by 
$\tilde d^{P+ \vert \bar Q\vert}$ and the choice of vertices of open simple self-intersections is bounded by $2H_\theta$. 
This observation, together with inequalities (5.7), (5.8), and (5.9) 
proves the estimate (3.4). 

It should be noted that the upper bound $\tilde d^{P+ \vert \bar Q\vert}$ as well as the estimate of $\Upsilon$ (3.5) 
is of the multiplicative form  $p$- and $q$-edges. 
This implies the validity of these bounds without any  relation with  the orientation of $p$- and $q$-edges. 
Lemma 3.1 is proved.


\subsection{ Proof of inequality (4.18)}

There is one-to-one correspondence between the set $\Theta_{2s}$ of Dyck paths of $2s$ steps $\theta_{2s}$ and
the set $\CT_s$ of plane rooted  trees $T_s$ of $s+1$ vertices. 
This correspondence is given by the chronological run over the tree $T_s$  \cite{St}. 
Given a vertex $\u$ of $T_s$, we determine in natural way its exit degree as a number of exits from this vertex
at the marked instants of time. 
Given $\theta$ and $T(\theta)$, we determine the height of the tree by relation  $H_T = H_{\theta}$.

\vs
\noindent  {\bf Lemma 5.2.} {\it Given a natural $d\ge 1$, let us consider a subset $\CT_s^{(\rho; d)}\subset \CT_s $
of trees such that the the root vertex $\rho$ is of  the exit degree $d$.  
Then 
$$
\vert \CT^{(\rho; d)}_s \vert \le 2\,  e^{-d\eta }\,   t_{s-1}\, ,
\eqno (5.10)
$$
where $\eta = \log (4/3)$ and $t_s= \vert \CT_s \vert = \vert \Theta_{2s}\vert = (2s)! /(s! (s+1)!)$ is the Catalan number. }

\vs

Let is present the main lines of  the proof of this statement \cite{K} based on the recurrent relation 
$$
t_s = \sum_{j=0}^{s-1}\,  t_{s-1-j}\, t_j, \quad t_0=1.
\eqno (5.11)
$$
We denote by $\t^{(d)}_s = \vert \CT_s^{(\rho; d)}\vert$.
Then 
$$
\t^{(d)}_s =  \ \ \sum_{\stackrel{l_1,\dots,l_d\ge 0} {l_1+\dots+l_d= s-d} }\ \ 
t_{l_1}\cdots t_{l_d}.
\eqno (5.12)
$$
One can interpret (5.12) as the sum over sub-trees constructed on $d$ different roots. 

We rewrite  (5.12)  in the form
\begin{eqnarray*}
\t^{(d)}_s &=& 
\sum_{ v=0}^{s-d} \ \ 
\sum_{\stackrel{l_1,\dots,l_{d-2}\ge 0 }{l_1+\dots+l_{d-2}= s-d-v}}
\ \ t_{l_1} \cdots t_{l_{d-2}} \ \ 
\left(\  \sum_{\stackrel{j_1,j_2\ge 0}{j_1+j_2=v} } \ 
t_{j_1} t_{j_2}\right) 
\\
&=& 
\ \  
\sum_{\stackrel{l_1,\dots,l_{d-1}\ge 0}{l_1+\dots+l_{d-1}= s-d+1}}
\ \ t_{l_1} \cdots t_{l_{d-1}} \  -
\sum_{\stackrel{l_1,\dots,l_{d-2}\ge 0}{l_1+\dots+l_{d-2}= s-d+1}}
\ \ t_{l_1} \cdots t_{l_{d-2}},
\end{eqnarray*}
where we have used (5.11).
Then we obtain the following recurrence 
$$
\t_{s}^{(d)} = \t_{s}^{(d-1)} - \t_{s-1}^{(d-2)}, 
\eqno (5.13)
$$
with $d\ge 2,   \t_{s}^{(1)} = t_{s-1}$  
and $ \t^{(0)}_s =0$.
\vss

Relation (5.13) implies that 
$\check t_{s}^{(d)}\le \check t_{s}^{(2)}$. 
It follows from (5.11) that $\t_{s}^{(2)}= t_{s-1}$. Then 
$
\t_{s}^{(d)}\le t_{s-1}.
$
Rewriting (5.12) in the form 
$$
\t_{s}^{(d)} = \sum_{v=0}^{s-d} \t_{s-v-1}^{(d-1)} \, t_v,
$$
and taking into account that $\t_{s-v-1}^{(d-1)}\le t_{s-v-2}$, we find that for $d-2\ge 1$
$$
\t_{s}^{(d)}\le \sum_{v=0}^{s-m} t_{s-v-2} \, t_v
\le 
t_{s-1} - t_{s-2}.
$$
Using the elementary bound $t_{s-2} \ge  t_{s-1}/4$, we can write that 
$\t_{s}^{(d)}\le {3\over 4} \, t_{s-1}$. Using this inequality  and repeating 
the computations presented above,
we see that 
$\t_{s}^{(d)}\le \left( {3\over 4}\right)^2\, t_{s-1}$. 
By recurrence, we get the bound 
$$
\t_{s}^{(d)} \le \left({3\over 4}\right)^{d-2} \, t_{s-1}.
$$
This inequality implies (5.10).  Lemma 5.2 is proved. 

\vs
\noindent {\bf Corollary of Lemma 5.2.} {\it Let $\CT^{(d)}_s$, $d\ge 2$ be a collection of trees of $\CT_s$ such that there exists a vertex $\u_0$
with $\vert \Delta(\u_0)\vert \ge d$. Then
$$
\vert \CT^{(d)}_s\vert \le 2s \, e^{-d \eta}\,  t_s.
\eqno (5.14)
$$
}

\vs
{\it Proof.} Given  $l$, $0\le l\le s-d$, we consider the vertex $\u_0$ and construct in the upper half-plane
a tree of $s-l$ edges that have $d$ edges attached to $\u_0$, i.e. the element of the set $\CT^{(\u_0;d)}_{s-d-l}$.
Then we construct in the lower half-plane a tree of $l$ edges that has $\u_ˆ$ as the root vertex. 
Finally, we choose a new root vertex $\rho$ among the $l+1$ vertices of the lower tree and perform the chronological run starting from 
$\rho$. The Dyck path obtained correspond to an element from $\CT^{(d)}_s$ and all trees of this form can be obtained
by this procedure. Then 
$$
\vert \CT^{(d)}_s\vert = \sum_{l=0}^{s-d} \, (l+1)\, t_l \, \vert \CT^{(\u_0;d)}_{s-l}\vert  \le 2 s \, e^{-d\eta}\  \sum_{l=0}^{s-d} t_l \, t_{s-l-1}
$$
and (5.14) follows. 


\vs
\noindent {\bf Lemma 5.3.} { \it Let   $\CT^{(u)}_s = \{ T: \, T\in \CT_s, H_T=u\}$
and let 
$\check \CT^{(u;d)}_s$ be
a family of  trees $T\in \CT_s$ such that $H_T = u, T\in \check \CT^{(u;d)}_s$ and the root vertex of $T$ is of the exit degree $d$.
Then 
$$
\sum_{u=1}^s \, e^{ x u/\sqrt s}\  \vert \check \CT^{(u;d)}_s\vert \le 2d \, e^{- d\eta}\,  \sum_{u=1}^s \, e^{ x u/\sqrt s}\ 
 \vert  \CT^{(u)}_s\vert  =  2d \, e^{-d \eta }\, B_s(x),
 \eqno (5.15)
$$
where $B_s(x)$ is determined by relation (2.9). }

{\it Proof.}  Let us denote by $\hat t_s^{(u)}$ the number of plane rooted trees  $T_s$ such that $H_T\le u$.
We also denote by $\tilde t^{(u)}_s = \vert \CT_s^{(u)}\vert $
with obvious agreement that $\tilde t^{(u)}_s=0$ in the case when $u>s$. 
The set $\check \CT^{(u;d)}_s$ can be represented as the union of  $d$ disjoint subsets with respect on where the first
arrival at the height $u$  happens. Then  we can write that
 $$
 \vert \check \CT^{(u;d)}_s\vert = \sum_{j=1}^d\ \ \  \sum_{l_1 + \dots + l_d = s-d} \ 
 \ \hat t_{l_1}^{(u-2)}\cdots  \hat t_{l_{j-1}}^{(u-2)}\, \tilde t^{(u-1)}_{l_j} \, 
 \hat t_{l_{j+1}}^{(u-1)}\cdots \hat t^{(u-1)}_{l_d},
 $$
where the sum runs over all $l_i\ge 0$. Taking into account that $\hat t_l^{(u-2)} \le \hat t_l^{(u-1)}\le t_l$, 
it is easy to deduce from (5.13) inequality
$$
 \vert \check \CT^{(u;d)}_s\vert \le \  d \ \  \sum_{l_1 + \dots + l_d = s-d} \ 
 \ \tilde  t_{l_1}^{(u-1)}\,   t_{l_{2}}\,  t_{l_3} \cdots  t_{l_d}.
 $$
Then 
$$
\sum_{u=1}^s \, e^{ x u/\sqrt s}\  \vert \check \CT^{(u;d)}_s\vert \le
d\,  \sum_{u=1}^s \, e^{ x u/\sqrt s}\ \ 
 \sum_{l_1 + \dots + l_d = s-d} \ 
 \ \tilde  t_{l_1}^{(u-1)}\,   t_{l_{2}}\,  t_{l_3} \cdots  t_{l_d} 
 $$ 
$$
\le d\,  \sum_{L=0}^{s-d}\ \left({1\over t_L}   \sum_{u=1}^L e^{ x u/\sqrt a}\ \tilde t _L^{(u-1)} \right)\  t_L 
\sum_{l_2+\dots + l_d = s-d-a} \ t_{l_2} \cdots t_{l_d}
$$
$$
\le d B_s(x)\  \sum_{L=0}^{s-d} \ t_L\ \  \sum_{l_2+\dots + l_d = s-d-L} \ t_{l_2} \cdots t_{l_d}.
\eqno (5.16)
$$
\vs
\noindent Remembering  (5.10) and  (5.12), we deduce from   (5.16)  the bound (5.15).

\vs
\noindent {\bf Lemma 5.4.} {\it Inequality (4.18) is true.}

\vs
{\it Proof.} To prove Lemma 5.4, we have to generalize the result of Lemma 5.3 to the case when the vertex $\upsilon_0 $
such that  $\vert \Delta(\upsilon_0)\vert \ge d$ is not necessary the root vertex. We can construct the corresponding set
$\CT^{(u;d)}_s$ by the following procedure.

Let us take a number of edges, say $k$ and construct the line graph from the root vertex
$\rho$ to the vertex $\upsilon_0$. Clearly, the vertices of the obtained graph can be ordered by the chronological run over it. 
Then we consider  the $k+1$ vertices less or equal to $\u_0$ and  construct $k+1$ sub-trees on them using
 $a$ edges.
Using remaining $s-d-k-a$ edges,  we  construct
$d+k$ subtrees on $d$ vertices of $\Delta(\u_0)$ and $k$ vertices greater than $\u_0$. Certainly, the trees we construct
have to verify the condition $H_T=u$. 

\vs

Then we can split the obtained set $ \CT^{(u;d)}_s(k)$ into three disjoint subsets in dependence where
the first arrival at the height $u$ appears. We  write that 
$$
\vert \CT^{(u;d)}_s(k)\vert  = \vert \CT^{(u;d)}_s(k;1)\vert +\vert \CT^{(u;d)}_s(k;2)\vert +\vert \CT^{(u;d)}_s(k;3)\vert,
\eqno (5.17) 
$$
where we determine 
$$
\vert \CT^{(u;d)}_s(k;1) \vert = \sum_{j=0}^k \ 
\ \sum_{a,b,c} \ \   \ \sum_{l_0+\dots + l_k=a} \ \ \ 
\hat t_{l_0} ^{(u-1)}\cdots \hat t_{l_{j-1}}^{(u-j)} \, \tilde t_{l_j}^{(u-j)}\ \hat t_{l_{j+1}}^{(u-j-1)}\cdots \hat t_{l_{k+1}}^{(u-k-1)}
$$
$$
\times \ \ \sum_{l_1'+\dots+l'_d=b} \ \hat t_{l'_1}^{(u-k-1)}\cdots \hat t_{l'_d}^{(u-k-1)} \ \  
\ \sum_{l''_1 + \dots + l''_{k} = c} \ \hat t^{(u+1-k)}_{l''_1} \, \hat t^{(u+2-k)}_{l''_2} \cdots \hat t^{(u)}_{l''_k},
\eqno (5.18)
$$
where the sum over $a,b,c$ is such that $a+b+c=s-d$. Variables 
$\CT^{(u;d)}_s(k;2)$ and $\CT^{(u;d)}_s(k;2)$ are determined in obvious way by similar to (5.18) expressions. 
One can  estimate the sub-sum 
$$
S_1(k) = \sum_{u=1}^s e^{xu/\sqrt s } \ \vert \CT^{(u;d)}(k;1)\vert
$$
by repeating the reasoning used in the proof of Lemma 5.3. We take a particular value of $l_j=L$ in (5.18)
and  remove the hats and superscripts from all other factors $t$. 
Then we change the order of summation and use the bound
$$
{1\over t_L} \ \sum_{u=1}^{s} e^{xu/\sqrt s} \, \tilde t^{(u-j)}_l \le B_L(x)\le B_s(x).
$$
It is easy to see  that the sum 
$$
T^{(d)}_s(k) = \ \sum_{a,b,c} \ \   \ \sum_{l_0+\dots + l_k=a} \ \ \ 
 t_{l_0}\cdots    t_{l_{k+1}}
\ \ \sum_{l_1'+\dots+l'_d=b} \  t_{l'_1}\cdots  t_{l'_d} \ \  
\ \sum_{l''_1 + \dots + l''_{k} = c} \ t_{l''_1} \,   \cdots  t_{l''_k}
$$
represents the number of all  trees that have the vertex $\upsilon_0$ with $\vert \Delta(\u_0)\vert \ge d$
at the distance $k$ from the root vertex.
Then
$$
\sum_{k=0}^{s-d} S_1(k) \le \sum_{k=0}^{s-d} (k+1) B_s(x) 
T_s^{(d)}(k) \le s B_s(x) \sum_{k=0}^{s-d} T_s^{(d)}(k) \le 2s^2 \, e^{-d\eta} \, t_s,
\eqno (5.19)
$$ 
where we have used the bound (5.14).  The sums $\sum_k S_2(k)$ and $\sum_k S_3(k)$ can be estimated 
by the right-hand side of (5.19).
Lemma 5.4 is proved.

\subsection{BTS-instants and proper and imported cells}

In this subsection we study the properties of  walks $w_{2s}$ 
that give us tools to control the cases when the graphs $\Gamma(w_{2s})$
have vertices of large exit degree; we  follow mostly the 
lines of  \cite{KV}.


\begin{figure}[htbp]
\centerline{\includegraphics[width=12cm,height=4cm]{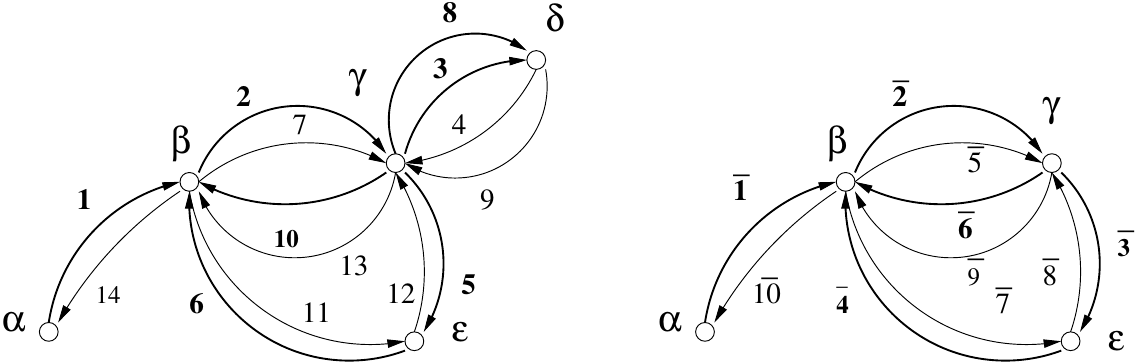}}
\caption{\footnotesize{ The graph   $\Gamma(w_{14})$  and its reduced counterpart
$ \bar \Gamma$}}
 \end{figure}

Given a walk $w_{2s}$, we determine  the following procedure of
reduction that we denote by $\cal P$: find an instant of time
$1\le t< 2s $ such that the step $(t-1,t)$ is marked and
$w_{2s}(t-1) = w_{2s}(t+1)$; if it exists, consider a new walk
$w'_{2s-2} = {\cal P}(w_{2s}) $ given  by a sequence
$$
w'_{2s} = (w_{2s}(0), w_{2s}(1), \dots , w_{2s}(t-1), w_{2s}(t+2), \dots, w_{2s}(2s)).
$$
 Repeating  this procedure as many times as it is
possible, we obtain the resulting walk denoted by $\bar W(w_{2s})$.

  Let us note that
$\bar W(w_{2s})$ is again a walk that can be transformed into the
new ordinary  one by renumbering  the values of $\bar W(t)$, $t\ge 1$ keeping their order.
Then we get the graph
 $\bar \Gamma(w_{2s})= \Gamma(\bar W(w_{2s}))$ that can be considered  as a sub-graph of $g(w_{2s})$,
$$
\CV(\bar \Gamma(w_{2s}))\subseteq \CV(\Gamma(w_{2s})),\quad
 \CE(\bar \Gamma( w_{2s}))\subseteq \CE(\Gamma(w_{2s}))
$$
and assume that the edges of $\bar \Gamma(w_{2s})$  are ordered
according to those  of $\Gamma (w_{2s})$.

 On Figure 2, the graph of the walk $w_{14}$ is shown as well as its reduced counterpart $\bar \Gamma (w_{14})$.

\vs 
\noindent {\bf Definition 5.1}. {\it Given a walk $w_{2s}$, we consider a
vertex  $\b\in \CV(\Gamma (w_{2s}))$ and refer to the marked
arrival edges  $(\a,\b)\in \CE(\Gamma (w_{2s}))$ as to  the  primary (or proper)
cells of $w_{2s}$ at $\b$. If $\b\in \CV(\Gamma (\bar W))$ with $\bar W
= \bar W(w_{2s})$, then we call the non-marked arrival edges
 $(a',\b)\in \CE(\Gamma (\bar W))$ the imported cells of $w_{2s}$ at $\b$. }

\vs As we will see later, in order to control the exit degree of a
vertex of a walk with typical $\theta$, one needs to take into
account the number of imported cells at this vertex. The main
observation here is that the presence of the imported cells is
closely related with the breaks of the tree structure performed by
the walk. 
Let us introduce the notion of the instant of broken tree structure.

\vs 
\noindent {\bf Definition 5.2.} {\it Any walk  $\bar W=\bar W(w_{2s})$
contains at least one instant $\bar \eta $ such that the  step
$(\bar \eta-1,\bar \eta)$ is marked and the step
 $(\bar \eta, \bar\eta+1)$ is not.
 We call such an instant $\bar  \eta$ the instant
of broken tree structure (or the BTS-instant of time) of the walk
$\bar W$. Passing back to the non-reduced walk $w_{2s}$, we
consider the edge $e(\eta)$ that corresponds to the edge
 $e(\bar\eta)\in \CE(\bar W)$ and refer to the instant $\eta$ as the
BTS-instant of the walk $w_{2s}$. }

\vs 
It is clear that   if the arrival instant $ \t$ is the
BTS-instant of the walk $w_{2s}$, then $\t$ is the open instant of
self-intersection of $w_{2s}$. The walk on Figure 2 contains two BTS-instants $\eta_1 = 6$
and $\eta_2=10$. These correspond to the 
instants of open self-intersections. 
The vertex $\gamma$ has one proper cell given by $t_1 = 2$ and one
imported cell $t_2= 7$. Regarding the Dyck path  $\theta(w_{14})$ and its tree $T_7$, 
we see that the edges $e(3)$ and $e(8)$ of $\Gamma (w_{14})$ represent different exit clusters of $T_7$. 
One can say that the edge  $e(8)$ is imported at $\gamma$ because in $T_7$ it does not belong
to the same exit cluster as the $e(3)$. This justifies the use of the term "imported cell". 

\vss
 Given  $\b\in \CV(w_{2s})$, we say that the 
BTS-instants $\eta_i$ such that \mbox{$w_{2s}(\eta_i)= \b$} are the
$\b$-local BTS-instants. All other BTS-instants are referred to as
the $\b$-remote BTS-instants.
 
\vss 
\noindent {\bf Lemma 5.5.}
 {\it Let us denote by  $ L^{(\setminus \b)}_{w_{2s}}$  the number
of all $\b$-remote  BTS-instants \mbox{of the walk $w_{2s}$}. Then
the number of all imported cells at $\b$  denoted by
$J_{w_{2s}}(\b)$ is bounded as follows}
$$
J_{w_{2s}}(\b) \le  L^{(\setminus \b)}_{w_{2s }} +
\kappa_{w_{2s}}(\b).
 \eqno (5.20)
$$

 {\it Proof.} Let us consider the reduced walk
 $\bar W = \bar W(w_{2s})$ and  $\b\in \CV(g(\bar W))$. We introduce the
function $\Phi_\b(t;\bar W)$ determined as the number of
$t$-open edges attached to $\b$;
$$
\Phi_\b(t;\bar W) = \# \{i: m(\a_i,\b;t) = 1 ({\hbox{mod }} 2)
\}
$$
and consider how this function changes its values at the instants
of time when $\bar W$ arrives at $\b$. The following
considerations show that this value can be changed by $ 0,+2$ and
$ -2$ only.

\vss
I. The value of $\Phi_\b$  stays unchanged  in the following two cases.

a) The first situation happens when  the walk $\hat W$ leaves $\b$
by a non-marked edge and arrives at $\b$ by a marked edge. Then
the corresponding cell is the primary one and we do not care about
it.

b) The second situation occurs when the walk $\hat W$ leaves $\b$
by a marked step $(x,x+1)$ and arrives at $\b$ by a non-marked
step $(y-1,y)$. Then the interval of time $[x+1,y-1]$ contains at
least one BTS-instant of time. This is the first instant $\tau$
when the non-marked step follows immediately after the marked one.
It is clear that $w_{2s}(\tau)\neq \b$ and therefore $\tau$ is the
$\b$-remote BTS-instant.

It should be noted that another such interval $[x'+1,y'-1]$
contains another $\eta'$ that obviously differs from $\eta$;
$\eta'\neq\eta$.  This is because  any couple of such time
intervals $[x+1,y-1]$, $[x'+1,y'-1]$ has an empty intersection.
Then each imported cell of the type (Ib) has at least one
corresponding $\b$-remote BTS-instant and the sets of the
BTS-instants that correspond to different intervals do not
intersect.

\vss
II. Let us consider the arrival instants at $\b$ when the value of
$\Phi_\b$ is changed.

a) The change by $+2$ takes place when the walk $\bar W$ leaves
$\b$ by a marked edge and arrives at $\b$ by a marked edge also.

b) The change by $-2$ occurs in the opposite case when $\hat W$
leaves $\b$ with the help of the non-marked edge and arrives at
$\b$ by a non-marked edge.

During the whole walk, these two different passages occur the same
number of times. This is because $\Phi_\b=0$ at the end of the
even closed walk $\bar W$. Taking into account that the number of
changes by $+2$ is bounded by the self-intersection degree
$\kappa_{\bar W}(\b)$, we conclude that the number of imported
cells of this kind  is not greater than
 $\kappa_{\bar W}(\b)$.

To complete the proof, we have to  pass back from the reduced walk
$\bar W(w_{2s})$ to the original $w_{2s}$. Since the number of
imported cells of $w_{2s}$ and the number of BTS-instants of
$w_{2s}$ are uniquely determined by $\bar W (w_{2s})$, and
$$
\kappa_{\bar W}(\b)\le \kappa_{w_ {2s}}(\b),
$$
then (5.20) follows for any $\b\in \CV(\Gamma(\bar W))$ as well as
to $\Gamma (w_{2s})$. 
If \mbox{$\b\notin \CV(\Gamma (\bar W))$,} then
$J_{w_{2s}}(\b)=0$ and (5.20) obviously holds.
 Lemma 5.5 is proved.

\vss \noindent 
{\bf Corollary of Lemma 5.5.} {\it Given a vertex $\b$ of the
graph of $w_{2s}$, the number of primary and imported cells
$G(\b)$ at $\b$ is bounded
$$
G(\b)\le 2\kappa(\b) +L,
 \eqno (5.21)
$$
 where $L$ is the
total number of the BTS-instants performed by the walk $w_{2s}$.}

\vs 
{\it Proof.} The number of the primary cells at $\b$ is given
by the $\kappa(\b)$. The number of the imported cells $J(\b)$ is
bounded by the sum $\kappa(\b) + L^{{(\setminus \b)}}$, where
$L^{(\setminus \b)}\le L$. Then (5.21) follows.

 \vss
Recently and independently from our work, 
the walks with BTS-instants are used in paper \cite{S1} and then in 
\cite{FS,S1,S2} as the base for a new combinatorial description 
of even closed walks. This approach is fairly powerful and promising. Also one should point out
a work \cite{M}, where the walks of broken tree structure were used in the study of 
traces of discrete Laplace operator. 


\subsection{Truncated random variables and proof of Theorem 2.1}

In the present subsection, we prove the following statement: 
 if random variables $\{ a_{ij}, 1\le i\le j <\infty\}$ are independent identically distributed random variables that 
have symmetric probability distribution such that 
$$
\E \vert a_{ij}\vert^{2(\phi+\d_0)}<\infty
\eqno (5.22)
$$
for some positive $\d_0<1$,  then 
 $$
 \P\left(  A^{(n)} \neq \hat A^{(n)} \ {\hbox{i. o.}} \right)=0,
 \eqno (5.23)
 $$ 
where $A^{(n)}$ is determined by (2.1) and the real symmetric matrix $A^{(n)}$ is such that 
$(\hat A^{(n)})_{ij} = {1\over \sqrt n} \hat a^{(n)}_{ij}, \ 1\le i\le j\le n$ with 
\mbox{$\hat a^{(n)}_{ij} = a_{ij} I_{(-U_n,U_n)}( a_{ij})$}
and 
$$
U_n= n^{1/\phi - \d}\  \ \ {\hbox{for any given }}\ \   0< \delta\le { \d_0\over  \phi(\phi +1)}.
\eqno (5.24)
$$

To prove (5.23), we use  the slight modification of the standard arguments of probability theory \cite{BY}. 
Condition (5.22) implies convergence of  the series 
$$
\sum_{k=1}^\infty P\{ \vert a_{ij}\vert^{2\phi +2\d_0} \ge k\}.
$$
{Taking} this into account and denoting $k= 2^{2m}+u$, we can write that
$$
\sum_{m=1}^\infty  \sum_{u=1}^{3 \cdot 2^{2m}-1}  \P\left( \vert a_{ij}\vert ^{2\phi +2\d_0} \ge 2^{2m}+u\right)\ge 
3\sum_{m=1}^\infty \  2^{2m}   \P\left( \vert a_{ij}\vert ^{2\phi +2\d_0} \ge 2^{2m+2}\right)
$$
$$
\ge 3 \sum_{m=1}^\infty \  2^{2m} \, \P\left( \vert a_{ij}\vert \ge 4\cdot 2^{m/(\phi +\d_0)}\right)
$$
and conclude that the last series converges.
From the other hand, we have 
$$
\P\left( A^{(n)}\neq \hat A^{(n)} \right) \le 
\sum_{m=k}^\infty \P\left( \cup_{2^{m-1}\le n\le 2^m} \ \cup_{1\le i\le j\le n} \ 
\left\{ \omega: \, \vert a_{ij}\vert \ge n^{1/\phi- \delta} \right\}\right)
$$
$$
\le \sum_{m=k}^\infty \P\left(  \cup_{1\le i\le j\le 2^m} \ 
\left\{ \omega: \, \vert a_{ij}\vert \ge 2^{m/\phi- m \d }\right\}\right)
$$
$$
\le \sum_{m=k}^\infty
2^{2m} \, \P\left( \vert a_{ij}\vert \ge 2^{m/\phi -m \d}\right)
$$
that vanishes as $k\to\infty$ provided $\phi \d< \d_0/(\phi+\d_0)$. 
Then (5.23) holds.

\vss
In fact, to prove relation (4.36), we do not need so strong estimates as above.
It is clear that 
$$
\P (\overline  {\cal O} _n) \le \sum_{1\le i\le j\le n}  \P\left\{ \vert a_{ij}\vert > n^{1/6-\d}\right\} \le
n^{2} { \E \vert a_{ij}\vert ^{12 +2\d_0}\over (n^{1/6-\d})^{12+2\d_0}}.
$$
Then (4.36) follows under condition that  $0<\d<\max\{1/6, \d_0/43\}$ with $\d_0\le 1$.

\section{ Summary}

In present work, we have studied the high moments of large Wigner random matrices with symmetrically distributed
entries  with the help of the improved version of the method proposed by Ya. Sinai and A. Soshnikov \cite{SS1,SS2,S}.
This improvement consists of two main counterparts: first, in paper \cite{KV} the Sinai-Soshnikov method 
has been completed by the rigorous study of the classes of walks with large number of instants of broken tree 
structure; second, in the spirit of paper  \cite{K}, the multiple edges were considered as the layers put one the basic structure
of the first-passage graph.

\vss
The new approach developed allows one to consider Wigner random matrices $A^{(n)}$ whose entries are given 
by independent identically distributed random variables that have a finite number of moments. 
In particular, we have proved the existence and universality of limits of the moments 
of random matrices with truncated elements, 
$\hat M_{2s_n}^{(n)} = \E \T (\hat A^{(n)})^{2s_n}$ when    $n$ and $s_n$ infinitely increase 
such that  $s_n = \chi n^{2/3}$ with $\chi>0$.

This can be considered as a step toward the proof of the universality of the Tracy-Widom law
in the Wigner ensemble of random matrices whose elements have a number of moments finite.

\vs 

The sufficient condition $V_{12+2\d_0}<\infty$ is imposed essentially  because of  the rate $s_n=\chi n^{2/3}$
that dictates the use of the truncation constant $U_n = n^{1/6-\d}$  with appropriate choose of $\d$. 
Regarding the main estimates of the paper,  one can easily observe that if $s_n = O(\log n)$, it is sufficient
to use the truncation variable of the form  $U_n= n^{1/2-\d}$ and the main  technical estimate  (2.8) remains true.
This result, together with the observations of subsection 5.7,  shows that the condition $V_{4+2\d_0}<\infty$ is sufficient 
for the convergence $\Vert A^{(n)}\Vert\to 2v$  to hold with probability 1 as $n\to\infty$. This condition is close to the 
optimal one $V_4<\infty$ determined in \cite{BY}. 
This means that our method could be also used to detect the optimal conditions for the convergence of
$\hat M_{2s_n}^{(n)}$ in the asymptotic regime $s_n = \chi n^{2/3}$.  

\vs It is argued in \cite{S} that the study of the correlation functions of the variables $L_{2s}^{(n)}$ (1.2) could be reduced
to the study of the moments $\E L_{2s'}^{(n)}$ with properly chosen $s'$. Then the universal upper bound of the form (2.7) 
would imply the existence of the universal limit of the corresponding correlation functions. Therefore, one can expect that
by using the technique developed in the present paper with possible use of additional arguments, 
one can broad the validity of the universality conjecture
confirmed in \cite{S} for the eigenvalue distribution of the Wigner ensemble whose entries satisfy condition (1.4).

{\bf Acknowledgements.} {The financial support of the research grant ANR-08-BLAN-0311-11 "Grandes Matrices Al\'eatoires" is gratefully acknowledged.}


\begin{thebibliography}{99}

\bibitem{BY} Z.D Bai and Y. Q. Yin, Necessary and sufficient conditions for the almost sure convergence of the largest eigenvalue of Wigner matrices, \emph{Ann. Probab.}  {\bf 16}  (1988), 1729-1741.

\bibitem{Baik} 
J. Baik, P. Deift, and K.  Johansson, On the distribution of the length of the longest increasing subsequence of random permutations, 
\emph{J. Amer. Math. Soc.} {\bf  12} (1999), 1119--1178.

\bibitem{Bol} B. Bollob\'as, 
\emph{Random graphs}, Academic Press, London, 1985. 

\bibitem{Boug} 
Ph.  Bougerol and Th. Jeulin, Paths in Weyl chambers and random matrices,
\emph{Probab. Theory Related Fields} {\bf  124} (2002),  517--543.

\bibitem{B} E. Br\'ezin et al. 
\emph{Applications of Random Matrices in Physics,} NATO Sciences Series II: Mathematics, Physics
and Chemistry, Vol. 221, Springer, Berlin, 2006.





\bibitem{FS} O. N. Feldheim and A. Sodin, A universality result for the smallest eigenvalues of certain sample 
covariance matrices, 
\emph{ Geom. Funct. Anal.}  {\bf 20}  (2010),   88Ð123.


\bibitem{FK} Z. F\"uredi and J. Koml\'os, The eigenvalues of random symmetric matrices, 
\emph{Combinatorica}, {\bf 1} (1981), 233-241.



\bibitem{G} V. Girko, \emph{Spectral properties  of random matrices}, Nauka, Moscow (1988) (in Russian)

\bibitem{J} 
K.  Johansson, Non-intersecting paths, random tilings and random matrices. 
\emph{Probab. Theory Related Fields} {\bf  123} (2002),  225--280.

\bibitem{Ok} 
A. Okounkov,  Random matrices and random permutations,
\emph{Internat. Math. Res. Notices}  (2000), 1043--1095.

\bibitem{K} A. Khorunzhy, Sparse random matrices: spectral edge and
statistics of rooted trees, \emph{Adv. Appl. Probab.} {\bf 33} (2001), 124-140.


\bibitem{KKPS} A. Khorunzhy, B. Khoruzhenko, L. Pastur, and M. Shcherbina, The large-n limit in statistical mechanics
and the spectral theory of disordered systems,  in: \emph{Phase Transitions and Critical Phenomena} {\bf Vol. 15},
pp.  74-239, Academic Press, London, 1992.
 

\bibitem{KM} O. Khorunzhiy and J.-F. Marckert, Uniform bounds for exponential moment of maximum 
of a Dyck path, \emph{Electr. Commun. Probab.}  {\bf 14}   (2009),  327--333.

\bibitem{KSV} O. Khorunzhy, M. Shcherbina, and  V. Vengerovsky, Eigenvalue distribution of large weighted random graphs,  
\emph{J. Math. Phys.}  {\bf 45}  (2004),  1648--1672.

\bibitem{KV} O. Khorunzhiy and V. Vengerovsky, Even walks and estimates of high moments 
of large Wigner random matrices, preprint, 2008.

\bibitem{M}  M. L. Mehta, \emph{Random Matrices}, 
Third edition. Pure and Applied Mathematics (Amsterdam), 142. Elsevier/Academic Press, Amsterdam, 2004, 688 pp.


\bibitem{Mn} P. Mn\"ev, Discrete path integral approach to the Selberg trace formula for regular graphs,
\emph{Commun. Math. Phys.} {\bf 274} (2007), 233-241.


\bibitem{P} L. Pastur, On the spectrum of random matrices, \emph{ Theor. Mathem. Physics} {\bf 10} (1972)

\bibitem{R} A. Ruzmaikina, Universality of the edge distribution of the
eigenvalues of Wigner random matrices with polynomially decaying distributions of
entries, \emph{Commun. Math. Phys.} {\bf 261} (2006), 277-296.


\bibitem{SS1} Ya. Sinai and A. Soshnikov, Central limit theorem for traces of large symmetric matrices with independent matrix elements, \emph{Bol. Soc. Brazil. Mat.} {\bf 29} (1998), 1-24.

\bibitem {SS2} Ya. Sinai and A. Soshnikov, A refinement of Wigner's semicircle law in a neighborhood
of the spectrum edge for random symmetric matrices, \emph{Func. Anal. Appl.} {\bf 32} (1998), 114-131.



\bibitem{S1} A. Sodin, Random matrices, non-backtracking walks, and the orthogonal polynomials,
{ \it J. Math. Phys.} {\bf 48} (2007), 123503.

\bibitem{S2} A. Sodin, The Tracy-Widom law for some sparse random matrices,
\emph{Journal of Statistical Physics} {\bf 136} (2009), 834-841.

\bibitem{S3} A. Sodin, The spectral edge of some random band matrices,
\emph{Ann. Math.} {\bf 172} (2010), 2223-2251



\bibitem{S} A. Soshnikov, Universality at the edge of the spectrum in Wigner random matrices,
\emph{Commun. Math. Phys.} {\bf 207} (1999), 697-733.

\bibitem{St} R. P.  Stanley, \emph{Enumerative Combinatorics},
Vol. 2. Cambridge University Press, 1999.



\bibitem{TW} C.Tracy and H. Widom, On orthogonal and symplectic matrix ensembles,
\emph{Commun. Math. Phys.} {\bf 177} (1996), 727-754.

\bibitem{W} E. Wigner, Characteristic vectors of bordered matrices with infinite dimensions,
\emph{Ann. Math.} {\bf 62} (1955), 548-564


\end{thebibliography}
\end{document}